\documentclass[11pt,leqno]{amsart}
\usepackage{framed}
\usepackage{a4wide}
\usepackage{amssymb}
\usepackage{amsmath}
\usepackage{amsthm}
\usepackage{mathrsfs}
\usepackage{euscript}
\usepackage[all]{xypic}
\usepackage[usenames]{color}
\usepackage{xr-hyper}
\usepackage{hyperref}

\usepackage{comment}
\usepackage{color}

\definecolor{shadecolor}{gray}{0.875}
\specialcomment{personal}{\begin{shaded}}{\end{shaded}}

\usepackage[usenames,dvipsnames]{xcolor} 
\CompileMatrices

   1

 3

\newcommand{\kll}{[\![}
\newcommand{\krr}{]\!]}

\newcommand{\Rep}{\operatorname{Rep}}
\newcommand{\id}{\operatorname{id}}

\newcommand{\fp}{\ifmmode {\mathbb{F}_p}\else$\mathbb{F}_p$\ \fi}
\newcommand{\zp}{\ifmmode {\mathbb{Z}_p}\else$\mathbb{Z}_p$\ \fi}
\newcommand{\z}{\mathbb{Z}}
\newcommand{\zpMod}{\ifmmode\mbox{$\zp$-Mod}\else$\zp$-Mod \fi}
\newcommand{\Mod}{\ifmmode\mbox{$\Lambda$-Mod}\else$\Lambda$-Mod \fi}
\newcommand{\La}{\ifmmode\Lambda\else$\Lambda$\fi}

\newcommand{\Hom}{{\mathrm{Hom}}}

\newcommand{\TLT}{T_{\pi}}

\renewcommand{\H}{\mathrm{H}}

\newcommand{\m}{\ifmmode {\frak m}\else$\frak m$ \fi}
\newcommand{\p}{\ifmmode {\frak p}\else${\frak p}$\ \fi}
\renewcommand{\P}{\ifmmode {\frak P}\else${\frak P}$\ \fi}
\newcommand{\e}{\ifmmode {\cal E}\else$\cal E$ \fi}
\newcommand{\G}{\ifmmode {\cal G}\else${\cal G}$\ \fi}
\renewcommand{\d}{\ifmmode {\cal D}\else${\cal D}$\ \fi}
\newcommand{\A}{\ifmmode {\cal A}\else${\cal A}$\ \fi}

\renewcommand{\projlim}[1] {{\lim\limits_{\stackrel{\displaystyle\longleftarrow}{#1}}}}

\renewcommand{\in}{\ \epsilon\ }
\newcommand{\Qp}{\ifmmode {{\Bbb Q}_p}\else${\Bbb Q}_p$\ \fi}
\newcommand{\qp}{\ifmmode {{\Bbb Q}_p}\else${\Bbb Q}_p$\ \fi}
\newcommand{\Q}{\ifmmode {\Bbb Q}\else${\Bbb Q}$\ \fi}
\newcommand{\Gal}{\mbox{\mbox{\rm Gal}}}

 \newcommand{\CC}{\mathbb{C}}

 \newcommand{\ZZ}{\mathbb{Z}}

 \renewcommand{\cL}{\mathcal{L}}
 
 \newcommand{\cO}{\mathcal{O}}

\renewcommand{\cR}{\mathcal{R}}

 \newcommand{\Zp}{\ZZ_p}

 \newcommand{\Cp}{\CC_p}
\newcommand{\rA}{\mathscr{A}}

\newcommand{\mf}{\mathfrak}

\newtheorem*{conj}{Conjecture}
\newtheorem{thm}{Theorem}
\newtheorem{cor}[thm]{Corollary}

\theoremstyle{definition}

\newtheorem{theorem}{Theorem}[section]
\newtheorem{corollary}[theorem]{Corollary}

\newtheorem{example}[theorem]{Example}

\usepackage{graphics}

\usepackage[latin1]{inputenc}

\usepackage{times}
\usepackage[T1]{fontenc}

\title[Explicit Reciprocity Laws] 
{Explicit Reciprocity Laws in Iwasawa Theory\\{\rm - A survey with some focus on the Lubin-Tate setting -}}

\author[Otmar Venjakob] 
{Otmar Venjakob}
\address{Universit\"{a}t Heidelberg,  Mathematisches Institut,  Im Neuenheimer Feld 288,  69120
Heidelberg,  Germany,
 http://www.mathi.uni-heidelberg.de/$\,\tilde{}\,$venjakob/}
\email{venjakob@mathi.uni-heidelberg.de}


\date{\today}%

\begin{document}

\begin{abstract}
Starting from Gau{\ss}' and Legendre's quadratic reciprocity law we want to sketch how it gave rise to the development of higher and generalized reciprocity laws and over all explicit reciprocity formulas in Iwasawa theory.
\end{abstract}

\maketitle

\section{Introduction}
The whole article is aimed to serve as an introduction to the different meanings of {\it explicit reciprocity laws or formulas}. Therefore we confine ourselves here to just give a short guideline through it. In section \ref{sec:gauss} we recall Gau{\ss}' reciprocity\footnote{We prefer the German spelling ''Gau{\ss}''  of this personal name, although non-German-speakers will be more familiar with ''Gauss'' perhaps.} law and link it to the quadratic Hilbert symbol. We then proceed in section \ref{sec:higherHilbert} to introduce higher Hilbert symbols in the context of Galois cohomology together with variants in characteristic $p$ or for formal groups. In section \ref{sec:explicitformulas} we explain some explicit formulas which calculate some of these higher Hilbert symbols. While already here towers of local fields play a crucial role, we only in section \ref{sec:PR} take the full Iwasawa theoretic point of view introduced by Perrin-Riou, Kato (and Coleman), by which we mean the {\it compatibility of cup-products in Galois cohomology with de Rham duality in $p$-adic Hodge theory}. As preparation and transition we recall in section \ref{sec:BK} the (dual) exponential map of Bloch and Kato and their version of an explicit reciprocity law for the representation $\Qp(r),$ for $r\geq1.$ Finally, in Iwasawa theory it has become common practise to also call statements like  \[Log_{BK/PR}(res_p(EulerSystem))=\mathcal{L}_p\]   expressing an $p$-adic $L$-functions $ \mathcal{L}_p$ as image of some global Euler system under the restriction map $res_p $ to local Galois/Iwasawa cohomology followed by the Bloch-Kato/Perri-Riou (dual) exponential/Big logarithm map $ Log_{BK/PR}$ as {\it explicit reciprocity laws} as we will explain in section \ref{sec:ES}\footnote{David Loeffler pointed out to me that this arguably    inaccurate terminology was probably first used in the  survey paper \cite{BCDDPR}: they refer to  Kato's work in which he {\it uses} a reciprocity law for formal groups (in the above mentioned traditional sense of the term) to deduce a formula for the image under the dual Bloch-Kato exponential map of his Euler system class; intentionally or accidentally they may have (mis)applied the name of the intermediate result as the name of the final theorem. Massimo Bertolini responded to me: {\it My recollection (somewhat vague after 30 years) is that the non-traditional terminology of calling {\it explicit reciprocity laws} also the logarithmic images of Euler systems emerged in the early 90's, in connection with Kato's work. Kato used an explicit reciprocity law (in the original sense) for higher dimensional local fields to prove his formula connecting special values of Hasse-Weil $L$-functions to the logarithms of his classes. After a while our community started to call this type of formula also an {\it explicit reciprocity law}. As Henri Darmon points out, one can find a similar pattern in the paper by Coates-Wiles \cite{coates-wil77}, where they used an {\it explicit reciprocity law} to establish their main theorem. By taking a quick look at their paper, I have the impression that a terminological distinction between the {\it explicit reciprocity law} and the result on special value of $L$-functions (theorem 29) was maintained.}
In any case the link between both, reciprocity laws and special $L$-values, is so strong - as we try to convince the reader in section \ref{sec:ES} - that this usage seems fully justified - over all when combined with David Loeffler's observation in subsection \ref{sec:lneqp} below.}. In the last section \ref{sec:eps} we quickly mention the impact of explicit reciprocity laws towards the  $\varepsilon$-isomorphism conjecture.\\
While we try to look at the general picture we sometimes focus on recent developments regarding the Lubin-Tate setting.

The survey is written mainly from the perspective of reciprocity laws. Nevertheless we aim at pointing out how interwoven the historical development has been with the aspect of special $L$-values and $p$-adic $L$-functions. Thus another option for reading this paper is to start with (or switch directly after section \ref{sec:explicitformulas} to) section \ref{sec:ES}.

Needless to say that the presentation follows an increasing gradient of abstraction and complexity according to the topics being discussed starting from elementary number theory and ending up to current research in this field.

We found the existing surveys
\cite{vostokov} on the history and \cite{li} on Kato's explicit reciprocity very helpful when preparing this manuscript which arose from different colloquium style talks the author had given on this subject.

{\bf Acknowledgements:}  I am very grateful to Denis Benois and Laurent Berger for their valuable comments and for giving advice concerning the existing literature. To David Loeffler, Henri Darmon and Massimo Bertolini I would like to express my gratitude for  discussions about the historic development of the terminology {\it explicit reciprocity law}. I thank   Muhammad Manji for allowing to contribute Example \ref{Manji}. I also would like to express my gratitude to Marlon Kocher, Rustam Steingart and Max Witzelsperger for helping to  improve the exposition and the correction of many typos. We   acknowledge funding by the Deutsche Forschungsgemeinschaft (DFG, German Research Foundation)  under TRR 326 {\it Geometry and
Arithmetic of Uniformized Structures}, project number 444845124.

%

%

%

\section{{ Gau{\ss}'} Reciprocity Law and the quadratic Hilbert symbol}\label{sec:gauss}

 In the preface of the wonderful book \cite{lemmer} we found the following to the point characterization of the {\it Quadratic Reciprocity Law}, which we recall immediately in subsection \ref{sec:GRL} below: \\

{\sc Erich Hecke (1923)\cite[p.59]{hecke1923}:}\\

\begin{quote}
  Modern number theory dates from the discovery of  {\bf the reciprocity law.} By its form it still belongs to the theory of rational numbers, as it can be formulated entirely as a simple relation between rational numbers; however its content points beyond the domain of rational numbers.\\
\end{quote}

{\sc Emma Lehmer (1978)\cite[p.\ 467]{lehmer}:}\\

\begin{quote}
... that the famous Legendre law of quadratic reciprocity, of which over 150 proofs are in print, has been generalized over the years [...] to the extent that it has become virtually unrecognizable.
\end{quote}

 The aim of this article is to reconfirm the reader with regard to the first quotation about the importance of  Legendre's and Gauss' quadratic reciprocity for the history of algebraic number theory. Most readers will immediately agree with the second statement - one motivation to write this article was also for the author to better understand the link between this quadratic reciprocity law and Perrin Riou's (generalised) reciprocity law. Nevertheless we shall try to uncover the common thread behind the developments which lead mathematicians to call those generalisations again {\it explicit reciprocity law } or {\it formula}.

\subsection{Legendre symbol}

The diophantine equation
\[X^2+pY=a\]
for $a\in \mathbb{Z}$ and an odd prime $p$ with $(p,a)=1$ has a solution in $\mathbb{Z}^2$ if and only if
\[X^2=\overline{a}\in \mathbb{F}_p^\times\]
has a solution in $\mathbb{F}_p$, i.e.,\ if $ a$ is a square there. This lead Legendre to define the following symbol which later has been extended by Jacobi allowing also non-prime natural numbers in the "denominator":



%

\begin{align*}
  \mbox{{\sc Legendre/Jacobi}-Symbol}  \phantom{mmm}\left(\frac{a}{p}\right):&=\left\{
                                                              \begin{array}{ll}
                                                                1, & \hbox{$\overline{a}\in(\mathbb{F}_p^\times)^2;$} \\
                                                                -1, & \hbox{otherwise.}
                                                              \end{array}
                                                            \right.
\end{align*}
Another description is given by {\sc Euler's}-criterion
\[\left(\frac{a}{p}\right)\equiv a^{\frac{p-1}{2}}\mod p.\]
Indeed, since  $\mathbb{F}_p^\times $ is cyclic of order  $p-1$   the following sequence is exact:
\begin{equation}\label{f:legendre}
\xymatrix@C=0.5cm{
  0 \ar[r] & (\mathbb{F}_p^\times)^2 \ar[rr]^{ } && \mathbb{F}_p^\times \ar[rr]^{\frac{p-1}{2}} && \{1,-1\} \ar[r] & 0, }
\end{equation}
where we consider
$\mu_2=\{-1,1\}$ as subgroup of $\mathbb{F}_p^\times.$

\subsection{Gau{\ss}' reciprocity law}\label{sec:GRL}

\phantom{mmmmmmmmmmmmmmmmmmmmmm}


The following quadratic reciprocity law had already been discovered by {\sc Euler} in 1744 and was  formulated by {\sc Legendre} in 1788. It was Gau{\ss} who first presented a complete proof of it:\\

\begin{minipage}[t]{\textwidth}
{\bf Reciprocity Law I ({\sc Gau{\ss}} 1801) :}
Let $l\neq p$ be odd prime numbers. Then we have
\[ \left(\frac{l}{p}\right)=(-1)^{\frac{l-1}{2}\frac{p-1}{2}}\left(\frac{p}{l}\right).\]
{\bf Supplement:}{  $\left(\frac{-1}{p}\right)=(-1)^{\frac{p-1}{2}},\;\; \left(\frac{2}{p}\right)=(-1)^{\frac{p^2-1}{8}}.$}\\

 {\bf Slogan:} \\
{\it $l$ is a square modulo $p$ if and only if $p$ is a square   modulo $l$  -- unless $l\equiv p\equiv 3\mod 4.$}\\

\end{minipage}
%

%

Even more general the Reciprocity Law holds for odd, pairwise coprime  natural numbers $m,n$ instead of $l,p$ using the Jacobi-symbol. This formulation of the reciprocity law\footnote{See \cite[Thm.\ 2.28]{lemmer} for variants.} explains in an obvious way its name, the two "fractions" on the left and right hand side being reciprocal to each other. This literal meaning of "reciprocity" gets lost in all generalisations of it we are going to discuss later on, which is one reason why it is hard to recognise the link to its origin.

Using the multiplicativity of the Legendre symbol in the upper variable and setting  $p^*:=(-1)^{\frac{p-1}{2}}p$ one obtains the following {\bf equivalent formulation:}
\begin{equation}\label{f:equiFormulation}
  \left(\frac{p^*}{l}\right)=\left(\frac{l}{p}\right).
\end{equation}

As mentioned already there are more than 150 proofs in the literature - Gauss himself already had included at least   6   in his Disquisitiones Math. We would like to sketch at least one proof here, which already is based on certain principles of class field theory. According to the last formulation \eqref{f:equiFormulation} it is natural to consider the quadratic extension $\Q(\sqrt{p^*})$ of $\Q$ which is visibly an abelian extension of $\Q.$ By Kronecker's Jugendtraum all abelian extensions of $\Q$ are contained in a cyclotomic extension. By considerations of discriminants it turns out that $\Q(\sqrt{p^*})$ is contained in $\Q(\zeta_p),$ wehre $\zeta_p$ denotes a primitive $p$th root of unity. Now the group $\mathbb{F}_p^\times$ occurs not only as the units of the residue field $\mathbb{F}_p$ at the prime $p$, but also as the Galois group $G(\Q(\zeta_p)/\Q) $ of this cyclotomic extension by mapping $a\in \mathbb{F}_p^\times$ to $\sigma_a$ satisfying $\sigma_a(\zeta_p)=\zeta_p^a$. By definition $\sigma_{\bar{l}}$ is nothing else than the
 Frobenius automorphism $Fr_l$ at $l$.
This leads to the following argumentation: Firstly,  an easy calculation in quadratic extensions shows that $Fr_l(\sqrt{p^*})=\left(\frac{p^*}{l}\right)\sqrt{p^*}$. Secondly, interpreting   the defining sequence \eqref{f:legendre} of the Legendre symbol also in terms of Galois groups attached to the embedding $\Q(\sqrt{p^*})\subseteq \Q(\zeta_p)$ we have the following commutative diagram with exact rows:\footnote{$G(\Q(\sqrt{p^*})/\Q) $  is the unique quotient of the cyclic group $ G(\Q(\zeta_p)/\Q)$ of order $2.$}

\phantom{mmmmmmmmmmmmmmmmmmmmmmm} { $Fr_l\in$}
\[\xymatrix@C=0.5cm{0 \ar[r] & G(\Q(\zeta_p)/\Q(\sqrt{p^*})) \ar[d]^{\cong} \ar[r]^{} & G(\Q(\zeta_p)/\Q) \ar[d]^{\cong}\ar[r]^{ } & G(\Q(\sqrt{p^*})/\Q) \ar[d]^{\cong}\ar[r] & 0\\
  0 \ar[r] & (\mathbb{F}_p^\times)^2 \ar[r]^{  } & \mathbb{F}_p^\times \ar[r]^{\left(\frac{\cdot}{p}\right) } & \{1,-1\} \ar[r] & 0. }\]
\phantom{mmmmmmmmmmmmmmmmmmmmmmmmmm} { $\bar{l}\in$}

From these two facts we obtain the equivalences
\[\phantom{mm}\left(\frac{p^*}{l}\right)=1\Leftrightarrow (Fr_l)_{\mid \Q(\sqrt{p^*})}=\id \Leftrightarrow \bar{l} \in (\mathbb{F}_p^\times)^2 \Leftrightarrow \left(\frac{l}{p}\right)=1,\]
which imply the reciprocity law.
From the diagram it becomes clear that we can read $\left(\frac{-}{p}\right) $ (or equivalently $\left(\frac{p^*}{-}\right) $) as a homomorphism in $\Hom(G(\Q(\zeta_p)/\Q),\mu_2)$ and for the proof it was crucial to evaluate them at the Frobenius automorphisms. Let's keep this in mind for later generalisations.

\subsection{Hilbert symbol}

A central principle in number theory  is the {\bf Local-Global-Principle.} Consider the absolute values on  $\mathbb{Q}:$ Besides the (real) absolute value $|-|_\infty$ there is for each prime $p$ the $p$-adic norm $|-|_p$ defined as
$|p^m\frac{a}{b}|_p:=p^{-m},$ if $(p,ab)=1.$ The completions of $\Q$ with respect to these norms show quite different behaviour:\\

\begin{center}
\begin{tabular}{ll}
$\mathbb{Q}_\infty=(\Q,|-|_\infty)^\wedge=\mathbb{R} $ & $\Qp:=(\Q,|-|_p)^\wedge$\\ \\
$\mathbb{R}^{>0}\xrightarrow{\log} (\mathbb{R},+) $ & $\mathbb{Q}_p^\times \xrightarrow{\log} \Qp,$
\end{tabular}
\end{center}
\vspace{0.5cm}

e.g.\ the subset
  $\mathbb{Z}_p:=\{z\in \Qp\mid |z|_p\leq 1\}$ forms again a subring in the $p$-adic world. The elements
$v\in\{p\mid \mbox{prime}\}\cup \{\infty\}$ are called the places of $\Q,$ a place $v\neq \infty$ is called {\em finite.}

The {\bf Quadratic Hilbert symbol} is now defined {\it locally} for each place $v$, respectively each local field $\Q_v$:
\phantom{mmmmmmmmmm}
\begin{align*}
  \left(\frac{m,n}{v}\right) &:=\left\{
                                 \begin{array}{ll}
                                   1, & \hbox{if $mX^2+nY^2=Z^2$ has non-trivial solution in $\mathbb{Q}_v$;} \\
                                   -1, & \hbox{otherwise.}
                                 \end{array}
                               \right.
\end{align*}
It defines  a symmetric, non-degenerate pairing
\[\left(\frac{-,-}{v}\right): \mathbb{Q}_v^\times/(\mathbb{Q}_v^\times)^2\times \mathbb{Q}_v^\times/(\mathbb{Q}_v^\times)^2 \to \mu_2.\]
%

By {\sc Hensel's }Lemma  the solvability of $X^2=\overline{a}\in \mathbb{F}_p^\times  $ is equivalent to that of  $ X^2\equiv a \mod p^n,$ for all $n\geq 1,  $ and in turn to  that of  $ X^2=a  \in\Zp$. This  can be used to easily conclude the identity
\begin{equation}\label{f:Hilbert=Legendre}
    \left(\frac{p,q}{p}\right)= \left(\frac{q}{p}\right)
\end{equation}
for   distinct odd positive primes numbers $p,q,$ see \cite[Prop.\ 2.26]{lemmer} for details and compare with \cite[V. Thm.\ (3.6)]{Neu} for a different proof.

It is an important insight that in number theory for the full picture one has to consider all places $v$ simultaneously as they are linked among each other, e.g., we have the\\

\centerline{{\bf
Closeness relation:}\phantom{mmmmmmm}  $\prod_v |a|_v=1$ for all $a\in\Q$.}
\vspace*{0.45cm}

This is a rather trivial consequence of the unique prime factor decomposition of integers In contrast the following statement is highly non-trivial:\\

{\bf Hilbert's Reciprocity Law II:} For $m,n\in \mathbb{Q}^\times$ it holds:
\begin{equation}\label{f:Hilbertglobal}
  \prod_{v\; \mbox{\footnotesize all places}} \left(\frac{m,n}{v}\right)=1.
\end{equation}

Indeed, it turns out that the laws  I and II are equivalent, for example the implication ''II $\Rightarrow$ I'' can be derived as follows:
\begin{equation}\label{f:LawII}
  1=\prod_v \left(\frac{p,q}{v}\right)=\left(\frac{p,q}{2}\right)\left(\frac{p,q}{p}\right)\left(\frac{p,q}{q}\right)=(-1)^{\frac{q-1}{2}\frac{p-1}{2}}\left(\frac{q}{p}\right)\left(\frac{p}{q}\right).
\end{equation}
Here we used the positivity of $p,q$ to see that the contribution at $\infty$ is one, moreover it is easy to verify that there is no non-trivial contribution from finite places distinct from $2,p,q.$ The Hilbert symbol at $2$ gives rise to the sign factor while we had seen already in \eqref{f:Hilbert=Legendre} how the other two Hilbert symbols specialize to the Legendre symbols.
For more details and the opposite implication we refer the interested reader to \cite[Thm.\ 2.28]{lemmer}.

\section{Higher Hilbert symbols and variants}\label{sec:higherHilbert}

The above proof of the Quadratic Reciprocity Law suggests the importance to study the group $\Hom(G_\Q,\mu_n)$ for varying natural numbers $n\geq 2$, which leads directly into class field theory and higher symbols as we will see in this section. Note that for a  finite extension $K$ and a discrete $G_K$-module $A$ with trivial $G_K$-action the group $\Hom_{cts}(G_K,A)$ of continuous group homomorphisms coincides with the first continuous group cohomology $H^1(G_K,A)=:H^1(K,A)$. In section \ref{sec:gauss} we had by chance that $\mu_2\subseteq \Q$ and hence $\Hom(G_\Q,\mu_2)=H^1(\Q,\mu_2).$ For general $n$ one has either to assume that $\mu_n\subseteq K$ or to work with $H^1(K,A).$

\subsection{Artinian  Reciprocity and different concepts of generalisations}

Gau{\ss}' Reciprocity Law is the beginning of {\bf class field theory}, which classifies and describes all  abelian  extensions of  global or local fields. The Artinian  Reciprocity map can be considered as generalisation of it, for example in the case of (finite)  Galois  extensions of local fields $K/F$ it comes as a surjective homomorphism
\[(-,K/F): GL_1(F)=F^\times  \twoheadrightarrow G(K/F)^{ab}\] with kernel the norm subgroup $N_{K/F}K^\times\subseteq F^\times.$ Instead of discussing  the shape of the global Artin map in general we just mention that in its ideal theoretic version the global Artin symbol $((q),\Q(\zeta_p)/\Q)$ attached to the principal ideal $(q)$ for the extension $\Q(\zeta_p)/\Q $ equals the arithmetic Frobenius automorphism $\sigma_q$, which sends $\zeta_p$ to $\zeta_p^q,$ whence its restriction
$((q),\Q(\sqrt{p^*})/\Q)$  to the extension $\Q(\sqrt{p^*})/\Q $ equals the Legendre symbol $\left(\frac{p^*}{q}\right)$ in the sense that $((q),\Q(\sqrt{p^*})/\Q)\sqrt{p^*}= \left(\frac{p^*}{q}\right)\sqrt{p^*},$ see \cite[after Cor.\ 3.10, Prop.\ 2.21]{lemmer}. Then Gau{\ss}' Reciprocity Law can be seen to be equivalent to the decomposition law of prime ideals in this quadratic extension, see \cite[I.Theorem (10.3), Prop.\ (10.5) and discussion afterwards,VI. Theorem (7.3)]{Neu}.

A {\it non-abelian generalisation} of the Artin reciprocity is the {\bf local (classical or $p$-adic) Langlands program}  relating roughly speaking certain $n$-dimensional representations of the absolute Galois group $G_F$ of $F$ to certain representations of the reductive algebraic groups $GL_n.$ For a recent survey on reciprocity laws and Galois representations, see \cite{wein}.

 {\it In this article}   we focus rather on {\bf generalisations of the quadratic Hilbert symbol}. To this aim we fix  a natural number $n\geq 2$, a finite extension
$ F/\Qp$ containing $\mu_n$, and we denote by $L/F$ the maximal abelian extension of exponent $n$ of $F$, i.e., such that $nG(L/F)=1.$ The {\it Kummer sequence}
\begin{equation*}\footnotesize\xymatrix@C=0.5cm{
  1 \ar[r] & {\mu_n} \ar[rr]^{ } && {\bar{F}}^\times \ar[rr]^{{\cdot}^n} && {\bar{F}}^\times \ar[r] & 1 }
\end{equation*}
induces by going over to the attached long exact sequence in Galois cohomology a canonical isomorphism
\[F^\times/(F^\times)^n\xrightarrow{\delta} H^1(F,\mu_n)\cong \Hom(G(L/F),\mu_n),\]
while the local Artin map induces an isomorphism
\[F^\times/(F^\times)^n\xrightarrow{(-,L/F)} G(L/F)\]
upon noting that the norm group $N_{L/F}L^\times$ of $L$ coincides with $(F^\times)^n.$
Then the {\bf $n$th Hilbert symbol} $ \left(\frac{\;,\;}{F}\right)_n $ is defined by the commutativity of the upper half of the following diagram, in which the first line is the Pontrjagin duality pairing, i.e., given by evaluating a group homomorphism from the right at an element from the left:
\begin{equation}\label{f:cup-product-Hilbert}\footnotesize
\xymatrix@C=0.5cm{
G(L/F)\ar@{}[r]|{\times} &{\Hom(G(L/F),\mu_n) }\ar[r]^-{{Pontrjagin} } & **[r]{\mu_n } \\
   F^\times/(F^\times)^n\ar[d]_{-\delta}\ar[u]_{(-,L/F)} \ar@{}[r]|{\times} &  F^\times/(F^\times)^n
   \ar[r]^{\left(\frac{\;,\;}{F}\right)_n}\ar[u]_{\delta}\ar[d]_{\delta} & **[r]{ \mu_n} \ar@{=}[u] \;\;\phantom{mm} \\
H^1(F,\mu_n)    \ar@{}[r]|{\times} &  H^1(F,\mu_n)\ar[r]^(0.3){\cup}   &
      H^2(F,\mu_n^{\otimes 2})\cong H^2(F,\mu_n)\otimes \mu_n \ar[u]_{\mathrm{inv}\otimes \mu_n}^{\cong}. }
\end{equation}
Unravelling this definition one obtains that
\[\left(\frac{a,b}{F}\right)_n:=\frac{(a,L/K)(\sqrt[n]{b})}{\sqrt[n]{b}}\]
for $a,b\in F^\times,$ see \cite[V Satz (3.1)]{Neu}. By \cite[Cor.\ (7.2.13)]{nsw}\footnote{f\"{u}r welche Normalisierungen der Klassenk\"{o}rpertheorie?} also the lower part of the diagram commutes, which establishes a relation between the Hilbert symbol and the cup product pairing for $\mathbb{G}_m$ respectively the Galois representation $\Qp(1)$ given by the cyclotomic character. Moreover, for $F=\Qp,$ $n=2$ we recover our quadratic Hilbert symbol  $\left(\frac{\;,\;}{\Qp}\right)_2=\left( \frac{\;,\;}{p} \right)$ as follows from \cite[V. Satz (3.6)]{Neu}.
In other words, we encounter two ways of generalisation of the quadratic Hilbert symbol:
\begin{enumerate}
\item[1.] The $n$th Hilbert symbol, which shall turn out to be related directly to the search of reciprocity laws concerning higher power residues, i.e., the  congruences $x^n\equiv a \mod p,$ with $n>2.$
\item[2.] Tate's local cup-product paring for different Galois representations. This topic will govern the Iwasawa theoretic development as discussed in section \ref{sec:PR}.
\end{enumerate}
For now we will introduce some variants of the higher Hilbert symbol.

\subsection{Schmid-Witt-Residue-Formula}\label{schmid-witt}
For a variant with
$\mathrm{char}(F)=p$ we consider the field $F=\mathbb{F}_p((Z))$. Replacing the Kummer sequence by the  {\it Artin-Schreier} sequence
\[\xymatrix@C=0.5cm{
  0 \ar[r] &  \mathbb{F}_p  \ar[rr]^{ } &&  {F}^{\mathrm{sep}} \ar[rr]^{\wp} &&  {F}^{\mathrm{sep}} \ar[r] & 0 }\]
with $\wp(x):=(Fr-\id)(x)=x^p-x$ induces a pairing
\[(\;,\;]: F^\times/(F^\times)^{p} \times F/\wp(F) \to \mathbb{F}_p\]
by the same recipe. It can be explicitly calculated by the
{\bf Schmid-Witt-Residue-Formula}
\[ (a,b]=\mathrm{Res}\left(b\frac{da}{a}\right)\]
with   $ \mathrm{Res}(
       (\sum_i a_i Z^i)dZ):= a_{-1}.$

 For $n>1$ there is similarly the Artin-Schreier-Witt pairing
 \begin{align*}
   [\;,\;)_F : W_n(F) \times F^\times & \longrightarrow  W_n(\mathbb{F}_p) ,
\end{align*}
 involving Witt-vectors $ W_n(F)$ of length $n,$  for which the Schmid-Witt formula gives again an explicit description, see \cite[\S 7]{SV15} for a version for {\it ramified} Witt-vectors, which will be crucial for the proof of Theorem \ref{thm:Kummermap} below. There is the classical article by Witt \cite{witt}. In the case of local fields of characteristic $p$, Parshin \cite{parshin} develops the higher class field theory using the generalized Witt pairing.\footnote{
For $n>1$ let $\mathscr{A}_L$ be  the $\pi_L$-adic completion of $o_L[[Z]][Z^{-1}],$ such that
  $\rA_L/\pi_L\rA_L=k((Z))=: K$,
   $ \Omega^1_{\mathscr{A}_L} = \mathscr{A}_L dZ$\\
    $d\log : o_L((Z))^\times \longrightarrow \Omega^1_{\rA_L}$ sends $f$ to $\frac{df}{f}.$
     Define the upper pairing in
\begin{equation*}
\xymatrix{
   W_n(\rA_L)_L \ar[d]_{\Phi_{n-1} }^{ }   \ar@{}[r]|{\times} & {o_L((Z))^\times} \ar[d]_{  d\log }\ar[r]^-{\{\ , \ \}} & o_L \ar@{=}[d]_{}^{} \\
     {\rA_L  } \ar@{}[r]|{\times} & {\Omega^1_{\rA_L}}\ar[r]^-{\mathrm{Res}} & o_L }
\end{equation*}
via the commutativity of the diagram. Moreover, there is a unique well defined bilinear pairing $(\ ,\ )$   such that the diagram
\begin{equation*}
\xymatrix{
   W_n(\rA_L)_L \ar@{->>}[d]_{W_n(\alpha_1)_L }   \ar@{}[r]|{\times} & {o_L((Z))^\times} \ar@{->>}[d]_{ \mod \pi_L }\ar[r]^-{\{\ , \ \}} & o_L \phantom{ } \ar@{->>}[d]^{\alpha_n } \\
   {W_n(K)_L} \ar@{}[r]|{\times} & {K^\times}\ar[r]^-{(\ ,\ )} & W_n(k)_L, }
\end{equation*}
is commutative.
\begin{thm}[Schmid-Witt-Formula for ramified Witt-vectors, {\cite[Thm.\ 7.16]{SV15}}]
The pairings $[\ ,\ )$ and $(\ ,\ )$ coincide.
\end{thm}}

\subsection{Lubin-Tate formal groups}

Another variant of the Hilbert symbol can be formulated for
(Lubin-Tate) formal groups $\mathcal{F}$  over a finite extension  $F$ of $\Qp$ attached to the prime  $\pi$ (and similarly for $p$-divisible groups). The sequence
\[\xymatrix@C=0.5cm{
  0 \ar[r] & \mathcal{F}[\pi^n] \ar[rr]^{ } && \mathcal{F}(\mathfrak{m}_{\bar{F}}) \ar[rr]^{[\pi^n]} &&  \mathcal{F}(\mathfrak{m}_{\bar{F}}) \ar[r] & 0 }\]
  induces a pairing
\[(\;,\;)_{\mathcal{F},n}: F^\times  \times \mathcal{F}(\mathfrak{m}_{{F}}) \to \mathcal{F}[\pi^n].\] We shall say more about the Lubin-Tate setting just before Theorem \ref{thm:wiles}.

So the culminating question is:\\

 {\bf How can one compute the Hilbert symbol $\left(\frac{\;,\;}{F}\right)_n$ (and its above variants) explicitly?}

\section{Explicit Formulas }\label{sec:explicitformulas}

\subsection{Tame case}

In the {\it tame case}, i.e., $p\nmid n$, the Hilbert symbol has a simple description:

Let $q$ be the cardinality of  $O_F/\pi_FO_F$ and consider the unique decomposition
\[ O_F^\times\cong \mu_{q-1}\times (1+\pi_FO_F),\;\; u\mapsto \omega(u)\langle u\rangle \]
with $ \omega(u)\in\mu_{q-1}$ and $\langle u\rangle \in 1+\pi_FO_F $;
let $v_F$ be the normalized valuation, i.e., $v_F(\pi_F)=1.$

For $a,b\in F^\times$ such that $\alpha=v_F(a)$, $\beta=v_F(b)$ it holds by \cite[V Satz (3.4)]{Neu}:
\[ \left(\frac{a,b}{F}\right)_n= \omega\left( (-1)^{\alpha\beta} \frac{b^\alpha}{a^\beta} \right)^{\frac{q-1}{n}}.\]

 In particular, for $a=\pi_F$ and $u\in O_F^\times$ we obtain that
 $\left(\frac{\pi_F,u}{F}\right)_n
=\omega(u)^{\frac{q-1}{n}}$ is the root of unity $\zeta\in\mu_n$ determined by $\zeta\equiv u^{\frac{q-1}{n}}\mod \pi_FO_F$, i.e., we have (see V.Satz (3.5) in (loc.\ cit.)):

\[\left(\frac{\pi_F,u}{F}\right)_n=1 \Longleftrightarrow \mbox{ $u$ is a $n$th power $\mod \pi_FO_F.$} \]

In other words the {\bf $n$th power residue symbol}
\begin{equation}\label{f:HigherPower}
\left(\frac{u}{\pi_F}  \right)_n:= \left(\frac{\pi_F,u}{F}\right)_n \mbox{ for $u\in  O_F^\times$ }
\end{equation}
  generalises the Legendre symbol $\left( \frac{a}{p}\right)$ which is reobtained by specializing to $F=\Qp$ and $n=2$. See VI \S 8 in \cite{Neu} for an extension to the definition of a global $n$th power residue symbol over arbitrary number fields $K$ (containing $\mu_n$) and the proof of a general reciprocity law (Theorem (8.3)) for it  as a consequence of the abstract formalism of class field theory, i.e., Artin reciprocity.  Indeed, one has again a product formula for $a,b\in K^\times$
\[\prod_v \left( \frac{a,b}{K_v}  \right)=1,\]
which results from the product formula $\prod_v(a,L_v/K_v)=1$ of (local) Artin maps for (principal id\`{e}les) $a\in K^\times.$
Similar as in \eqref{f:LawII} upon using \eqref{f:HigherPower} one obtains for $a,b\in K^\times$ prime to each other and to $n$ the general formula
\[\left(\frac{a}{b}\right)\left(\frac{b}{a}\right)^{-1}=\prod_{v\mid n\infty}\left( \frac{a,b}{K_v}  \right). \]

Theorem VI.(8.4) in (loc.\ cit.) or \eqref{f:LawII} explains how Gau{\ss}' reciprocity law is a special case. In other words: The long search in number theory for similar laws for $n$th power residue symbols lead to the discovery  of the Artin reciprocity law, which in turn eventually lead to a full explanation of the general reciprocity law of higher power residue symbols - except for the explicit calculation of the $n$th Hilbert symbol in the case of (wild) ramification, i.e., for $v\mid n,$ which we are going to discuss in the next subsection.

 \subsection{Wild case: Kummer, Artin-Hasse, Iwasawa, Wiles, Br\"{u}ckner, ...}
 Therefore henceforth we consider for odd $p$ the other extreme
\[n=p^k, \] i.e., we replace $n$ by $p^n$ and take $F:= K_n:=\mathbb{Q}_p(\zeta_{p^n})$ where  $\zeta_{p^n}$ primitive $p^n$th root of unity. Then
$\pi_n:=\zeta_{p^n}-1$ forms a prime element for the valuation ring of integers
$   o_{K_n}=\{|z|_p\leq 1\}$. We also introduce the trace map
 $Tr:=Tr_{K_n/\mathbb{Q}_p}:=\sum_{\sigma\in G(K_n/\mathbb{Q}_p)} \sigma$ for this extension $F/\Qp.$

It was Kummer \cite{kummer}, who found the following description of the Hilbert symbol for $n=1$, in which  differential logarithms of units show up: For  $\beta\in 1+\pi_1o_{K_1}$ he writes $\beta=f(\pi_1)$ for some $f\in\zp[[Z]]$ and performs calculus with regard to the variable $\pi_1.$ so, formally $d\log b$ amounts to the logarithmic derivative $d\log f(Z)=\frac{f'(Z)}{f(Z)}$ evaluated at $Z=\pi_1.$ In the same spirit the residuum $\mathrm{Res}_{\pi_1}$ will be taken with respect to $\pi_1.$

\begin{thm}[{\bf Kummer 1858}]\label{thm:kummer}
 For $\alpha, \beta\in 1+\pi_1o_{K_1}$ it holds:
\begin{align*}
  (\alpha,\beta)_{p}&=\zeta_{p}^{\mathrm{Res}_{\pi_1}(\frac{\log\alpha\cdot d\log\beta}{\pi_1^p})}= \zeta_{p}^{{\frac{1}{p}Tr(\zeta_p\log\alpha\cdot d\log\beta)} }
\end{align*}
\end{thm}

As Franz Lemmermeyer \cite{lemmer-overflow} points outs that "Kummer worked out the arithmetic of cyclotomic extensions guided by his desire to find the higher reciprocity laws; notions such as unique factorization into ideal numbers, the ideal class group, units, the Stickelberger relation, Hilbert 90, norm residues and Kummer extensions owe their existence to his work on reciprocity laws. His work on Fermat's Last Theorem is connected to the class number formula and the "plus" class number, and a meticulous investigation of units, in particular Kummer's Lemma, as well as the tools needed for proving it, his differential logarithms, which much later were generalized by Coates and Wiles." In his talk during the memorial conference of John Coates in Cambridge 2023 Andrew Wiles pointed out that indeed, Kummers work and Iwasawa's interpretation of it was crucial for both the development of his reciprocity law in the Lubin-Tate setting (Theorem \ref{thm:wiles} below) as well as the link to special-$L$ values (\eqref{f:CWder},\eqref{f:CWderelliptic}) and $p$-adic $L$-functions (\eqref{f:p-adicL}, \eqref{f:p-adicL-elliptic}), which we will discuss in section \ref{sec:ES}.

Only many years later Artin-Hasse \cite{artin-hasse} obtained a result for arbitrary $n.$

\begin{thm}[{\bf Artin-Hasse 1928}]
 For $ \beta\in 1+\pi_no_{K_n}$ it holds:
\begin{align*}
  (\zeta_{p^n},\beta)_{p^n}&=\zeta_{p^n}^{{\frac{1}{p^n}Tr(\log\beta)} } \\
  (\beta, \pi_{n})_{p^n}&=\zeta_{p^n}^{\frac{1}{p^n}Tr(\frac{\zeta_{p^n}}{\pi_n}\log\beta) }
\end{align*}
\end{thm}

Iwasawa \cite{iwasawa} contributed the following version, which specializes to Kummer's version for $n=1.$

\begin{thm}[{\bf Iwasawa 1968}]\label{thm:iwasawa}
   For $ \beta=(\beta_k)\in \projlim{k, Norm} K_k^\times$, $\phantom{}g_\beta\in \mathbb{Z}_p[[Z]]$ with $g_\beta(\pi_n)=\beta_n$ and   $\alpha\in 1+\pi_no_{K_n}$ we have:
\begin{align*}
  (\beta_n, \alpha)_{p^n}&=\zeta_{p^n}^{\frac{1}{p^n}Tr(\log(\alpha)\cdot D\log\beta) }
\end{align*}
with invariant  logarithmic derivation
\[D\log\beta=\left((1+Z)\frac{g_\beta'(Z)}{g_\beta(Z)}\right)_{\mid Z=\pi_n}.\]
\end{thm}
This result should also be compared with Coleman's more complete formula in \cite{col-di}.
Due to our interest towards Iwasawa theory we concentrate here on the above versions, but we would like to also mention:
Helmut Br\"{u}ckner \cite{brueckner} discovered in 1964 an explicit description by a residuum formula over not necessarily cyclotomic fields, see \cite[V. Theorem (3.7)]{Neu}. Another treatment has been given by Guy Henniart \cite{henniart}. See also the work of Sergei Vostokov et.\ al.\ \cite{vostokov, vostokov-explicit, vostokov-witt}. \\

Wiles \cite{wiles-explicit} generalized Iwasawa's result to the
{\bf Lubin-Tate setting:} we fix a  finite extension
$L/\mathbb{Q}_p$   with prime element  $\pi\in o_L$  as well as a Lubin-Tate formal group $\mathcal{F}=\mathcal{F}_\pi$ over $L$
   attached to $\pi$. We write  $[a]_\mathcal{F}(Z) \in o_L[[Z]]$, $a\in o_L,$ for the power series giving the $o_L$-action on it and $\log_\mathcal{F}$ for the logarithm of the formal group.
  The $\pi_L^n$-division points generate a tower of Galois extensions $L_n=L(LT[\pi_L^n])$ of $L$ the union of which we denote by $L_\infty$ with Galois group $\Gamma_L.$
   We let   $ T:=T_\pi\mathcal{F}=\projlim{n}\mathcal{F} [\pi^n]=o_L\eta$ denote its Tate-module with $o_L$-generator
$\eta=(\eta_n)$ and on which the  Galois action is described by the Lubin-Tate  character
 $\chi_{LT} :
 G_L \longrightarrow o_L^\times$. Furthermore, we write  $\partial_{\mathrm{inv}}$ for the invariant derivation with respect to $\mathcal{F}$ and $t_L:=\log_{\mathcal{F}}(Z)$, i.e., $\partial_\mathrm{inv}(f) = g_{\mathcal{F}}^{-1} f'$, where $g_{\mathcal{F}}$ is the formal derivative of  $\log_{\mathcal{F}}.$ As before let $Tr=Tr_{L_n/L}$ denote the corresponding trace map.
For  $L=\Qp,$ $\pi=p$, $\mathcal{F}=\hat{\mathbb{G}}_m$ this specializes to the above cyclotomic setting:\\

\begin{tabular}{lll}
{\bf Lubin-Tate} && {\bf cyclotomic}\\ \\
 $L/\mathbb{Q}_p$   &  degree $d$ extension of & {\color{blue}$\Qp$}\\
$o_L$ & integers & {\color{blue}$\Zp$}\\
$\pi_L\in o_L$& prime element & {\color{blue}$p$}\\
$q=|o_L/\pi_Lo_L|$ & cardinality of residue field & {\color{blue}$p=|\Zp/p\Zp|$}\\
 $\mathcal{F}=\mathcal{F}_\pi$ & Lubin-Tate formal group/$L$  &{\color{blue} $\hat{\mathbb{G}}_m$}\\
 & attached to $\pi=\pi_L$& \\
 $[a](Z) \in o_L[[Z]]$ & giving $o_L$-action, $a\in o_L$&{\color{blue}$(Z+1)^a-1$}\\
  $T_\pi\mathcal{F}=\projlim{n}\mathcal{F} [\pi^n]$ & Tate-module, $G_L$-action by & {\color{blue}$\Zp(1)$}, \\
   $\chi_{LT} :
 G_L \longrightarrow o_L^\times$ & Lubin-Tate  character & {\color{blue}$\chi_{cyc}$}\\
$\Omega$ & period of Cartier dual & \\
 $\varphi_L$ & $f(Z)\mapsto f([\pi_L](Z))$  & {\color{blue}$\varphi_\Qp$}\\
   $\psi_L$ & (almost)left inverse of $\varphi_L$ & {\color{blue}$\psi_\Qp$}
 \end{tabular}
\vspace*{0.5cm}

Then Wiles obtained the following generalisation.

\begin{thm}[{\bf Wiles 1978}]\label{thm:wiles}
 For $ \beta=(\beta_k)\in \projlim{k, Norm} L_k^\times$,\phantom{mmm} $\phantom{}g_\beta\in o_L[[Z]]$ with $g_\beta(\eta_n)=\beta_n$ and   $\alpha\in \mathcal{F}(\eta_no_{L_n})$ it holds:
\begin{align*}
  (\beta_n, \alpha)_{\mathcal{F},n}&=[{\frac{1}{\pi^n}Tr\left(\log_\mathcal{F}(\alpha) D\log g_\beta(\eta_n)\right) }]_\mathcal{F}(\eta_n)
\end{align*}
with   invariant logarithmic derivation $D\log g_\beta=\frac{1}{\log_\mathcal{F}'}\frac{g_\beta'}{g_\beta}.$
\end{thm}

While Coates and Wiles \cite{coates-wiles-explicit, coates-wil78, coates-wil77} were working on extending Iwasawa theoretic methods to CM-elliptic curves (as will be later discussed in section \ref{sec:ES}), Coleman systematized their methods. In particular, he discovered what one now calls Coleman power series, which have the property to interpolate a full norm-compatible system of local units for all layers of the Lubin-Tate tower simultaneously. We recall the injective ring endomorphism
\begin{align*}
  \varphi_L : o_L[[Z]] & \longrightarrow o_L[[Z]] \\
                     f(Z) & \longmapsto f([\pi_L](Z)) \ .
\end{align*}
Moreover, there is 
a unique multiplicative map $\mathscr{N} : o_L[[Z]] \longrightarrow o_L[[Z]]$ such that
\begin{equation*}
  \varphi_L \circ \mathscr{N} (f)(Z) = \prod_{a \in LT_1} f(a +_{\mathcal{F}} Z)  \qquad\text{for any $f \in o_L[[Z]]$}
\end{equation*}
(\cite{col1979} Thm.\ 11).

\begin{thm}[Coleman]\label{thm:Coleman}
For any norm-coherent sequence $u = (u_n)_n \in \varprojlim_n L_n^\times$ there is a unique Laurent series $g_{u,\eta} \in (o_L((Z))^\times)^{\mathscr{N} = 1}$ such that $g_{u,\eta}(\eta_n) = u_n$ for any $n \geq 1$. This defines a multiplicative isomorphism
\begin{align*}
  \varprojlim_n L_n^\times & \xrightarrow{\; \cong \;} (o_L((Z))^\times)^{\mathscr{N} = 1} \\
  u & \longmapsto g_{u,\eta} \ .
\end{align*}
\end{thm}

We now reformulate Iwasawa's theorem \ref{thm:iwasawa} in the way of \cite[\S 1.1]{kato-generalized}. Via the isomorphism $\z/p^m\z\cong \mu_{p^m},$ $\bar{1}\mapsto \zeta_{p^m},$ we interpret the $p^m$th Hilbert symbol as pairing
\[K_m^\times \times K_m^\times \to \z/\p^m\z.\]
Fixing $m\geq 0$ in the first argument and taking inverse limits in the second argument (with respect to norm maps) and in the target gives pairings
\begin{equation}\label{f:limitpairing}
 K_m^\times \times \varprojlim_n K_n^\times \to \zp, \;\;\;\mbox{ and }  (-,-):(K_m^\times)\otimes_{\zp}\Qp \times \varprojlim_n K_n^\times \to \Qp.
\end{equation}
Hence the non-degenerateness of the trace pairing implies the unique existence of a map $\lambda_m: \varprojlim_n K_n^\times   \to K_m$ making the following diagram commutative
\begin{equation}\label{f:Katodiag}
   \xymatrix{
    ( K_m^\times )\otimes_{\z}\Q\ar@{}[r]|{\times} & {\varprojlim_n K_n^\times \ar[d]_{\lambda_m}} \ar[r]^-{(-,-)} &{ \Qp } \ar[d]^{ } \\
     K_m \ar[u]_{\exp }\ar@{}[r]|{\times} & K_m \ar[r]^{Tr_{K_m/\Qp}} & {\Qp} ,  }
\end{equation}
where $\exp(a)$ is defined as $\exp(p^ka)\otimes p^{-k}$ for $k$ sufficiently large. More explicitly, $\lambda_m$ is the composite
\begin{equation}\label{f:lambdacomp}
  \varprojlim_n K_n^\times \xrightarrow{\mbox{Hilbert symbol}} \Hom_{cts}( K_m^\times, \Qp)\xrightarrow{\exp^*} \Hom_{cts}( K_m, \Qp)\xleftarrow[\cong]{Tr}K_m.
\end{equation}

Then the statement of Theorem \ref{thm:iwasawa} means that
\begin{equation}\label{f:lambdam}
  \lambda_m(\beta)=\frac{1}{p^m}(D\log g_\beta)_{\mid Z=\pi_n}.
\end{equation}
In \S 1.2 of (loc.\ cit.) Kato explains a relationship between cyclotomic units and some values of partial Riemann zeta functions, which gives a first hint that explicit reciprocity laws are closely related to Euler systems and special $L$-values or even $p$-adic $L$-functions - a connection we shall resume in section \ref{sec:ES}.
A similar reformulation is of course also possible for Wiles's theorem with $\lambda_m: \varprojlim_n L_n^\times \to L_m$ given by $\frac{1}{\pi_L^m}(D\log g_{-})_{\mid Z=\eta_n}.$ Note that $\lambda_m$ is nothing else than (part of) the Coates-Wiles homomorphism $\psi^1_{CW,m}$ in sections \ref{sec:BK} and \ref{sec:ES}.\\[0.5cm]

We close this section by giving comments to the literature and mentioning that there are lots of further formulas by various mathematicians - we are very grateful to Denis Benois for some guidance in this regard: \\

\textsc{ \textbf{Classical reciprocity laws:}}
''Classical'' results can be divided in two parts: Artin-Hasse type
formulas and Shavarevich-Br\"{u}ckner-Vostokov type formulas. The foundational papers for these approaches are the paper of Artin-Hasse \cite{artin-hasse} and
the paper of Shafarevich \cite{shafarevich}, respectively; see also Kneser \cite{kneser}. The link between these formulas can be understood in the framework of syntomic cohomology (see
below).

The following three important papers added  new techniques to the development:\\
  Sen's formula in \cite{sen}  is more general than the formulas of Artin-Hasse-
Iwasawa. More importantly, his approach uses the computation of the
Galois cohomology of the field $\mathbb{C}_p.$ In some sense, this is the first
application of methods of $p$-adic Hodge theory to explicit reciprocity laws.\\
In \cite{kato-explicit} Kato first uses syntomic cohomology to prove an explicit reciprocity law. In some sense, the explicit reciprocity law of Bloch-Kato (to be discussed in the next section) and other explicit reciprocity laws proved by
Kato take their origin in this paper. This paper gives a conceptual
explanation why Hilbert symbols can be computed using differential
forms. See also the appendix of \cite{gros} written by Kurihara in which  he shows that   both types of explicit formulas
for the Hilbert symbol (Artin-Hasse-Iwasawa type and Shafarevich-Br\"uckner-Vostokov type formulas) can be recovered by syntomic techniques.\\
Fontaine's appendix of \cite{perrin94} is the first paper where the relationship between the characteristic
$0$ and characteristic $p$ cases is explicitly established. For the further
development of this point, see the papers of Abrashkin  \cite{abrashkin97} and \cite{abrashkin97b}.\\

\textsc{ \textbf{Explicit reciprocity laws for formal groups:}}
The   paper \cite{tavares} by Tavares Ribeiro generalizes the result of Abrashkin (the condition that $\mu_{p^n}$  is contained in  the ground field is removed):
The approach   is different from Abrashkin's one and
introduces $(\varphi, \Gamma)$-modules in the false Tate curve extensions.

To sum up, the most general formulas for formal groups {\it over one-dimensional local fields} are obtained in the paper of Tavares Ribeiro (for Br\"uckner-Vostokov formulas, only formal groups over unramified local fields can be covered by this type of formulas) and Benois' paper \cite{benois-crelle} (for Artin-Hasse type formulas,
the ground field is arbitrary, but there are restrictions on the valuation of the second argument). Kolyvagin's formula in \cite{kolyvagin} is really explicit
only in the case of Lubin-Tate groups and coincides in that case with
Wiles' explicit reciprocity law. In the case of general formal group his formula depends on
some non explicit Galois invariants. These invariants can be explicitly
computed using $p$-adic periods accoding to \cite{benois-crelle}.\\
The paper of  Fl\'{o}rez \cite{flo} concerns Lubin-Tate formal groups {\it over higher
local fields.} The case of an arbitrary formal group in the higher dimensional case was studied by Fukaya \cite{fukaya}. Her result is a direct generalization
of \cite{benois-crelle}.  For the Hilbert symbol in the higher dimensional field case see also the article \cite{abrashkin-jenni} by
  Abrashkin and   Jenni.\\
e  In characteristic $p,$ concerning {\it  Formal Drinfeld Modules}, recently there has been remarkable progress by Eddine \cite{eddine2022norm,eddine2022explicit}. Her results are more general that the formulas of   Bars and
Longhi \cite{bars-longhi}. She follows the approach of Kolyvagin but applies it in the
characteristic $p$ case, where the theory of $p$-adic periods is not available.\\

\"{u}

We recommend  to consult the textbook \cite{fesenko-vostokov} as well as the surveys  \cite{vostokov}, \cite{bars-longhi} and the literature therein.

 \section{Bloch-Kato's reciprocity law}\label{sec:BK}

The explicit calculations of the higher symbols in section \ref{sec:explicitformulas} culminated into the interpretation \eqref{f:lambdam} based on the reformulation of the pairing into \eqref{f:lambdacomp}. This can also be rewritten as the statement that the class in $H^1(\Qp,\Qp(1))=\Hom_\Gamma(H_\Qp^{ab},\Qp(1))$ roughly associated with such $\lambda_m$ below is the image of $1$ under a transition map $\partial^1:\Qp\to H^1(\Qp,\Qp(1))$ arising from $p$-adic Hodge  theory. This is the starting point for Blochs and Katos variant.

Bloch and Kato \cite{bloch-kato} found a way to generalize the exponential map $\exp_\mathcal{F}$ of a formal group (as it shows up in Kato's reformulation \eqref{f:lambdam} of Iwasawa's and Wiles's explicit reciprocity law) to de Rham Galois representations. Built on this they were able to generalize the explicit reciprocity for $\Qp(1)$ (corresponding to the multiplicative formal group $\mathcal{F}=\hat{\mathbb{G}}_m$) to $\Qp(r).$

Let $K$ be a finite extension of $\Qp.$
For a continuous representation of $G_K$ on a finite dimensional $\mathbb{Q}_p$-vector space $V$ we write as usual\footnote{We refer the reader to \cite{FO} for the foundations in $p$-adic Hodge theory. }
\begin{align*}
   & D_{dR,K}(V) := (B_{dR} \otimes_{\mathbb{Q}_p}V)^{G_K} \supseteq D_{dR,K}^0(V):= (B_{dR}^+ \otimes_{\mathbb{Q}_p}V)^{G_K} \quad\text{and}  \\
   & D_{cris,K}(V) := (B_{max,\mathbb{Q}_p}\otimes_{\mathbb{Q}_p}V)^{G_K} \ .
\end{align*}
The quotient $tan_K(V) := D_{dR,K}(V)/D_{dR,K}^0(V)$ is called the tangent space of $V$.

Henceforth we assume that $V$ is de Rham. Then the usual Bloch-Kato exponential map $\exp_{K,V}:tan_K(V)\to H^1(K,V)$ can  be defined as follows. Apply the tensor functor $-\otimes_{\mathbb{Q}_p} V$ to the  exact sequence
\begin{equation}\label{f:FESQp}
    0 \to \mathbb{Q}_p \to B_{max,\mathbb{Q}_p}^{\phi_p = 1} \to B_{dR}/B^+_{dR} \to 0
\end{equation}
and take the (first) connecting homomorphism in the associated $G_K$-cohomology sequence.
Furthermore, the dual exponential map $exp^*_{K,V}$ is defined by the commutativity of the following diagram
\begin{equation}\label{f:dualexp}
\xymatrix{
  H^1(K,V) \ar[d]_{\cong} \ar[rrr]^{\exp^*_{K,V}} & & & D_{dR,K}^0(V) \ar[d]^{\cong} \\
  H^1(K,V^*(1))^* \ar[rrr]^-{(\exp_{K,V^*(1)})^*} & & & (D_{dR,K}(V^*(1))/D_{dR,K}^0(V^*(1)))^* ,  }
\end{equation}
where the left, resp.\ right, perpendicular isomorphism comes from local Tate duality, resp.\ from the perfect pairing
\begin{equation}\label{f:pairingDdR}
  D_{dR,K}(V) \times D_{dR,K}(\Hom_{\mathbb{Q}_p}(V,\mathbb{Q}_p(1)))  \longrightarrow D_{dR,K}(\mathbb{Q}_p(1)) \cong K ,
\end{equation}
in which the $D^0_{dR,K}$-subspaces are orthogonal to each other. Note that the isomorphism $K\cong D_{dR,K}(\mathbb{Q}_p(1))$ sends $a$ to $at^{-1}_{\mathbb{Q}_p} \otimes \eta^{cyc}$. Also, $(-)^*$ here means the $\mathbb{Q}_p$-dual, $t_{\mathbb{Q}_p}$ is Fontaine's period satisfying $g(t_{\mathbb{Q}_p})=\chi_{cyc}(g)t_{\mathbb{Q}_p}$ similarly as the basis $\eta^{cyc}$ of $\Qp(1).$

For a formal group $\mathcal{F}$ with $p$-adic Tate-module $T$ and $V:=T\otimes_\zp \Qp$ it is shown in \cite[Example 3.10.1]{BK} that the following diagram commutes

\begin{equation}\label{f:BKformal}
   \xymatrix{
     tan(\mathcal{F}) \ar[d]_{\cong } \ar[r]^{\exp_\mathcal{F}} & \mathcal{F}(o_K)\otimes_\z \Q \ar[d]^{Kummer} \\
     tan_K(V) \ar[r]^{\exp_{K,V}} & H^1(K,V).  }
\end{equation}
In particular, $\exp_{K, \Qp(1)}$ is induced by the classical $p$-adic exponential map $\exp$ while $\exp^*_{K, \Qp(1)}$ is the trivial map, because $tan_K(\Qp)=0.$\footnote{
Now assume that $V$ is in $Rep_L(G_K)$ and consider $K = L$ in the following. Let $\mathfrak{D}_{L/\mathbb{Q}_p} = \pi_L^s o_L$ be the different of $L/\Qp.$ Tensoring the $o_L$-linear  isomorphism
\begin{align}\label{f:different}
o_L & \xrightarrow{\;\cong\;} \Hom_{\mathbb{Z}_p}(o_L,\mathbb{Z}_p) \\
y & \longmapsto [x \mapsto \mathrm{Tr}_{L/\mathbb{Q}_p}(\pi_L^{-s}xy)]  \nonumber
\end{align} with $\mathbb{Q}_p$ gives the isomorphism of $L$-vector spaces
\begin{equation*}
  \tilde{\Xi}: L \cong \Hom_{\mathbb{Z}_p}(o_L,\mathbb{Z}_p)\otimes_{\mathbb{Z}_p}\mathbb{Q}_p \cong \Hom_{\mathbb{Q}_p}(L,\mathbb{Q}_p ) \ .
\end{equation*}
Since $\Hom_{\mathbb{Q}_p}(L,-)$ is right adjoint to scalar restriction from $L$ to $\mathbb{Q}_p$, and by using $\tilde{\Xi}^{-1}$ in the second step, we have a natural isomorphism
\begin{equation}\label{f:Ldual}
  \Hom_{\mathbb{Q}_p}(V,\mathbb{Q}_p) \cong \Hom_{L}(V,\Hom_{\mathbb{Q}_p}(L,\mathbb{Q}_p)) \cong \Hom_{L}(V,L) \ .
\end{equation}
Combined with \eqref{f:pairingDdR} we obtain the perfect pairing
\begin{equation}\label{f:pairingDdRoverL}
   D_{dR,L}(V) \times D_{dR,L}(\Hom_{L}(V,L(1))) \longrightarrow L
\end{equation}
with an analogous orthogonality property. Furthermore, similarly as in \cite[Prop.\ 5.7]{SV15} local Tate duality can be seen as a perfect pairing  of finite dimensional $L$-vector spaces
\begin{equation}\label{f:Tate-local}
  H^i(K,V) \times H^{2-i}(K,\Hom_{L}(V,L(1)) ) \longrightarrow H^2(K, L(1)) = L \ .
\end{equation}
Altogether we see that, for such a $V$, the dual Bloch-Kato exponential map can also be defined by an analogous diagram as \eqref{f:dualexp} involving the pairings \eqref{f:pairingDdRoverL} and \eqref{f:Tate-local} and in which $(-)^*$ means taking the $L$-dual.}

We can now interpret the map $\lambda_m$ in \eqref{f:lambdacomp} as
\begin{equation}\label{f:lambdacompBK}
  \varprojlim_n K_n^\times \xrightarrow{\mbox{Kummer map}} H^1_{Iw}(\zp(1))\xrightarrow{cor} \H^1( K_m, \Qp(1)) \xrightarrow{\exp^*_{BK,\Qp(1)}} D_{dR,K_m}^0(\Qp(1)) \xleftarrow{\cong}K_m
\end{equation}
and it is this shape, which will be suitable for generalizations. In particular, it should be  compared with the diagonal map in Corollary \ref{cor:KatoExplicit} below specialized to $\mathcal{F}=\hat{\mathbb{G}}_m$.

Recall from \eqref{f:lambdam} that the classical reciprocity law just states, that $\lambda_m$ is given by the first Coates-Wiles homomorphism. We will see later how the general Coates-Wiles homomorphisms arose in the study of $p$-adic $L$-functions in \ref{sec:p-adicL}. Bloch and Kato interpreted them as classes in $H^1(\Qp,\Qp(r))$ for $r\geq 1$ and reformulated for $r=1$ the classical reciprocity law by saying that this class arises as the image of $1$ under a certain transition map $\partial^r: \Qp \to H^1(\Qp,\Qp(r))$ arising from $p$-adic Hodge theory. We shall sketch here its generalisation in \cite[\S 8]{SV15} to the Lubin-Tate setting as introduced around Theorem \ref{thm:wiles}.


We define the Coates-Wiles homomorphisms in this context for $r \geq 1$ and $m \geq 0$ by\footnote{For $m>0$  one can extend the definition to $\varprojlim_n L_n^\times $ while for $m=0$ one cannot evaluate at $\eta_0=0$!}
\begin{align*}
   \psi_{CW,m}^r : \mathbb{U}(L_\infty):= \varprojlim_n o_{L_n}^\times & \longrightarrow L_m \\
                                   u & \longmapsto
    \frac{1}{r!\pi_L^{rm}}\left(\partial_{\mathrm{inv}}^{r-1}\frac{\partial_{\mathrm{inv}} g_{u,\eta}}{ g_{u,\eta}} \right)_{\mid Z=\eta_m} \ ,
\end{align*}
where $ g_{u,\eta}$ is the Coleman power series from Theorem \ref{thm:Coleman}.
Then the map
\begin{align*}
   \Psi_{CW,m}^r :  \mathbb{U}(L_\infty) & \longrightarrow L_mt_L^r \\
                                   u & \longmapsto \psi_{CW,m}^r(u)t_L^r
\end{align*}
is $G_L$-equivariant (it depends on the choice of $\eta$). In the following we abbreviate $\psi_{CW}^r := \psi_{CW,0}^r$ and $\Psi_{CW}^r := \Psi_{CW,0}^r$ . One might think about these maps in terms of the formal identity
\begin{equation}\label{f:CW1}
  \log g_{u,\eta} (Z) = \sum_r \psi_{CW}^r(u)t_L^r = \sum_r \Psi_{CW}^r(u)  \qquad\text{in $L[[t_L]] \subseteq B_{dR}$}
\end{equation}
or the identity
\begin{equation}\label{f:CW2} d\log g_{u,\eta}{(Z)}=\frac{dg_{u,\eta}{(Z)}}{g_{u,\eta}{(Z)}}= \sum_{r\geq 1}r\psi_{CW}^r(u) t_L^{r-1}dt_L \ .
\end{equation}

From class field theory we have an exact sequence
\[\xymatrix@C=0.5cm{
  0 \ar[r] & \mathbb{U}(L_\infty) \ar[rr]^{rec} && H_L^{ab} \ar[rr]^{ } && \hat{\mathbb{Z}} \ar[r] & 0 }\]
of $\Gamma_L$-modules, where the latter one has trivial action and $H_L^{ab}$ denotes the maximal abelian quotient of the absolute Galois group $H_L=G(\bar{L}/L_\infty)$ of $L_\infty$ with $\Gamma_L$-action induced by inner conjugation within $G_L.$ Therefore, for $r\geq 1$, one has isomorphisms \footnote{ Note that $H^i(\Gamma_L,L(\chi_{LT}^r))=0$ for $i,r\geq 1 $!}
\begin{align*}
  H^1(L,L(\chi_{LT}^r)) & \cong  H^1(L_\infty,L(\chi_{LT}^r))^{\Gamma_L}\\
    & =\Hom_{\Gamma_L}(H_L^{ab}, L(\chi_{LT}^r) ) \\
    & =\Hom_{\Gamma_L}( \mathbb{U}(L_\infty),L(\chi_{LT}^r) ).
\end{align*}
Since $\Psi_{CW}^r$ belongs to the latter group, we can interpret it as class in $ H^1(L,L(\chi_{LT}^r))  $. On the other hand the exact sequence
\begin{equation}\label{f:exseqBK}
    0  \longrightarrow Lt_L^r \longrightarrow Fil^r B_{max,L}^+ \xrightarrow{\pi_L^{-r}\phi_q-1} B_{max,L}^+ \longrightarrow 0
\end{equation}
from \cite[Lem.\ 8.2 (i)]{SV23}  induces, for any $r \geq 1$, the connecting homomorphism in continuous Galois cohomology
\begin{equation*}
  L = (B_{max,L}^+)^{G_L} \xrightarrow{\;\partial^r\;} H^1(L, L t_L^r)
\end{equation*}
and one can show \cite[Prop.\ 8.5]{SV15} that upon identifying $L t_L^r\cong L(\chi_{LT}^r)$ we have commutative diagram
\[\xymatrix{
  L=tan_L(L(\chi_{LT}^r)) \ar[rr]^{1-\pi^{-r}} \ar[dr]_{\exp_{L,L(\chi_{LT}^r) }}
                &  &    L \ar[dl]^{\partial^r}    \\
                & H^1(L, L(\chi_{LT}^r)) .                }\]

\begin{thm}[Bloch-Kato explicit reciprocity law {\cite[Thm.\ 2.1]{BK},\cite[Thm.\ 8.6 ]{SV15}}]\label{thm:BK}\footnote{Strictly speaking their explicit reciprocity law is Theorem 2.6 in (loc.\ cit.), which calculates the cup-product
\[H^1(K,\z/p^n\z(r))\times H^1(K(\mu_{p^n}),\z/p^n\z(1))\to H^2(K(\mu_{p^n}),\z/p^n\z(r+1))\] via differential forms using Coleman power series and
  which immediately implies the result stated here.}
For all  $a \in L$ and $r \geq 1$, we have the identities
\begin{equation*}
  {\partial}^r(a)= ar \Psi_{CW}^r
\end{equation*}
and
\begin{align*}
  \exp_{L,L(\chi_{LT}^r) }(a)& = (1-\pi_L^{{ -r}})ar \Psi_{CW}^r. \\
\end{align*}
interpreted as maps on $\mathbb{U}(L_\infty).$
\end{thm}

The following immediate consequence \cite[Cor.\ 8.7]{SV15} is again derived by considering duals and setting  $\mathbf{d}_r := t_L^rt_{\mathbb{Q}_p}^{-1} \otimes (\eta^{\otimes -r} \otimes \eta^{cyc})$, where $\eta^{cyc}$ is a generator of the cyclotomic Tate module $\mathbb{Z}_p(1)$, and $t_{\mathbb{Q}_p} := \log_{\mathbb{G}_m}([\iota (\eta^{cyc})+1]-1)$.

\begin{cor}[A special case of Kato's explicit reciprocity law] \label{cor:KatoExplicit}
For $r\geq 1$ the diagram
\begin{equation*}
  \xymatrix{
  {\varprojlim_n  o_{L_n}^\times} \otimes_{\mathbb{Z}} T^{\otimes -r} \ar[d]_{-\kappa\otimes \id} \ar[rrdd]^{\qquad ''(1-\pi_L^{-r})r\psi_{CW}^r({_-}) \mathbf{d}_r\, ''} &   \\
  H^1_{Iw}(L_\infty/L,T^{\otimes -r}(1)) \ar[d]_{\mathrm{cores}}  &   \\
  H^1(L, T^{\otimes -r}(1)) \ar[rr]^-{\exp^*} & & D^0_{dR,L}(V^{\otimes -r}(1)) = L \mathbf{d}_r ,  }
\end{equation*}
commutes, i.e., the diagonal map sends $u\otimes a \eta^{\otimes -r}$ to
\begin{equation*}
  a(1-\pi_L^{-r})r\psi_{CW}^r(u) \mathbf{d}_r = a \frac{1-\pi_L^{-r}}{(r-1)!} \partial_{\mathrm{inv}}^r \log g_{u,\eta}(Z)_{\mid Z=0} \mathbf{d}_r \ .
\end{equation*}
\end{cor}

We recommend also to compare with \cite{dS} and to consult \cite{saikia} for a survey on Bloch's and Kato's explicit reciprocity law. In that whole volume {\it The Bloch-Kato Conjecture for the Riemann Zeta Function} also the meaning for the Tamagawa number conjecture for the motive $\z(r)$ is discussed in detail.
The general case of Kato's explicit reciprocity law in this setting can be found in \cite{kato-lnm}, it has been generalized by Tsuji in \cite{tsuji}. As Laurent Berger pointed out his result is  one of rare ones which are not restricted to the $L$-analytic setting in the sense of subsection \ref{sec:abstract} below and it would be most desirable (but difficult) generalizing Theorem 5.3 in (loc.\ cit.)  to higher dimensional representations.

\section{Perrin-Riou's Reciprocity Law}\label{sec:PR}

 Perrin-Riou's Reciprocity Law stands for an explicit computation of

  \begin{center}
   {\it Tate's local  cup-product pairing}\\

     or\\

     {\it Iwasawa-cohomology pairing} \\[0,3cm]
 \end{center}

by means of the Big  Exponential map  $\Omega_{V^*(1)}$ (and/or the regulator map  $\mathcal{L}_V$)
 for a crystalline\footnote{E.g.\ representations attached to abelian varietes over $\Qp$ with {\it good} reduction. We refer the reader to \cite{FO} for the foundations in $p$-adic Hodge theory. } representation $V$ of   $G_{\qp}$
 \begin{center}\it
 in terms of easier to handle dualities from $p$-adic Hodge-theory.
 \end{center}
So this is a generalization of the interpretation of the Hilbert symbol as cup-product as in  \eqref{f:cup-product-Hilbert}, which justifies the use of ''reciprocity law''.

\subsection{Iwasawa cohomology and Big Exponential map}

 In order to explain Perrin Riou's point of view we first have to introduce
{\bf Iwasawa-cohomology} attached to the Galois tower of number fields
\[K \subseteq K_1 \subseteq \cdots \subseteq K_n\subseteq \cdots \subseteq K_\infty=\bigcup K_n\]
   with Galois group $\Gamma:=G(K_\infty/K)$ a $p$-adic Lie Group. We write
  $\Lambda(\Gamma)=\zp[[\Gamma]]$ for the Iwasawa algebra, i.e.,   completed group algebra. For a
 Galois stable lattice $T\subseteq V$ we define its Iwasawa-cohomology groups as

 \[H_{Iw}^i(T):=\projlim{n} H^i(K_n,T)\cong H^i(K,\Lambda(\Gamma)\otimes_{\zp} T)\] which are  finitely generated modules under $\Lambda(\Gamma)$ and important protagonists in the Iwasawa Main Conjectures together with its rational versions
 $H_{Iw}^i(V):=H_{Iw}^i(T) \otimes_{\zp} \Qp,$ see \cite{nek}, \cite{fukaya-kato}.

In \cite{perrin92,perrin94,perrin99} Perrin-Riou introduced the (big) exponential map for a cristalline representation $V$ of $G_\Qp$ satisfying $Fil^{-h}D_{cris}(V)=D_{cris}(V)$ for some integer $h\geq 1$ (such $h$ always exists):\footnote{Strictly speaking $\Omega_{V,h} $ is only defined on some submodule $ \bigg(D(\Gamma,\qp)\otimes_\qp D_{cris}(V)\bigg)^{\Delta=0}$, see e.g. \cite[Def.\ II.12]{berger-exp}. Moreover, in general, its image is only well-defined in $D(\Gamma,\qp)\otimes_{\Lambda} H_{Iw}^1(V)/V^{H_\Qp}. $ Up to replacing $V$ by an appropriate twist $V(j)$ one can always avoid these restrictions! }
\[\Omega_{V,h}: D(\Gamma,\qp)\otimes_\qp D_{cris}(V)\to D(\Gamma,\qp)\otimes_{\Lambda} H_{Iw}^1(V).\]
Here $D(\Gamma,\qp)$ denotes the distribution algebra of locally ($\Qp$)-analytic distributions of $\Gamma$ with coefficients in $\Qp.$ Note that all the twists $V(j)$ can be obtained as Galois equivariant quotients of $D(\Gamma,\Qp)\otimes_{\Qp} V.$
The big exponential map is uniquely characterized by the fact that it interpolates the Bloch-Kato exponential maps, the construction of which is recalled in section \ref{sec:BK},
\[\exp_{K_n,V(j)}:D^{K_n}_{cris}(V(j))\to H^1(K_n,V(j))\]
for all cyclotomic twists $V(j)$ of $V$, $j\in\mathbb{Z},$ and all integers $n\geq 0$ and satisfying a twist-behaviour with respect to the cyclotomic character.\footnote{$Tw_{V(j)}\circ\ \Omega_{V(j),h}\circ D=-\Omega_{V(j+1),h+1}$ in \cite[Thm.\ 3.2.3]{perrin94} or \cite[Rem.\ II.15]{berger-exp} without the sign.} Later Berger \cite[Thm:\ II.13]{berger-exp} found an easier direct construction based on use of $(\varphi,\Gamma)$- and Wach-modules. Up to tensoring $\Omega_{V,h}$ with the field of fractions of $D(\Gamma,\Qp)$ one can extend the definition to all $h\in\mathbb{Z}.$

\subsection{Different formulations of Perrin Riou's explicit reciprocity law}

In section \ref{sec:higherHilbert} we saw that one can interpret reciprocity laws via duality. This principle can be lifted to the level Iwasawa cohomology extended as modules over $D(\Gamma,\Qp)$:

Perrin-Riou formulates her explicit reciprocity law \cite[Conj.\ 3.6.4]{perrin94} as the compatibility of   Iwasawa duality with the (base change of the) crystalline (or de Rham) duality via her big exponential map as stated in (the lower part of)  diagram \eqref{f:explRec} below. As indicated in the upper part of the diagram, the Iwasawa duality pairing\footnote{This can either be defined by a limit process of Tate duality running up the tower $K_n$ as in \cite{perrin94} or by a residue pairing as in \cite[\S 4.2]{KPX} or \cite[\S 4.5.4]{SV23}.} $\{,\}_{Iw} $ interpolates the cup-product $<,>_{Tate}$ of usual Galois cohomology \`{a} la Tate. In this sense, one can use the crystalline pairing $ [,]_{D_{cris}}$ in order to calculate explicitly both the Iwasawa pairing and then the cup-product, hence her name {\it explicit reciprocity law} is in full accordance with the previous use in section \ref{sec:higherHilbert}. Perrin-Riou was able to prove this conjecture for $\Qp(r)$, $r\in\mathbb{Z}$, more generally it was shown  independently by {\sc Colmez 1998 \cite{colmez1998}, Benois 1998 \cite{benois}, Kato-Kurihara-Tsuji (unpublished)}

\begin{thm}[First formulation]\label{thm:PRI}
Let $V$ be crystalline. Then, for all integers $h, j$  we have a commutative diagram
\begin{equation}\label{f:explRec}\footnotesize
\xymatrix@C=0.5cm{
H^1(K_n,V^*(1-j) )\ar@{}[r]|{\times} &H^1(K_n,V( j))  \ar[r]^(0.45){<,>_{Tate}} & H^2(K_n,\mathbb{Q}_p(1))\cong\Qp \\
   {D(\Gamma,\qp)\otimes_{\Lambda} H_{Iw}^1(V^*(1)) }\ar[u]_{pr_n} \ar@{}[r]|{\times} &  {D(\Gamma,\qp)\otimes_{\Lambda} H_{Iw}^1(V)}
   \ar[r]^-{\{,\}_{Iw}}\ar[u]_{pr_n} & D(\Gamma,\qp) \ar@{=}[d]\ar[u]_{(-1)^j \int{\chi_{trivial}-}} \\
D(\Gamma,\qp)\otimes_{\qp} D_{cris}(V^*(1)) \ar[u]_{\Omega_{V^*(1),1-h}}  \ar@{}[r]|{\times} &  D(\Gamma,\qp)\otimes_\qp D_{cris}(V)\ar[r]^-{[,]_{D_{cris}} } \ar[u]_{\Omega_{V,h}} &
      D(\Gamma,\qp).}
\end{equation}
\end{thm}

Colmez \cite[Thm.\ 9]{colmez1998} proves the explicit reciprocity law in a different, but equivalent version, our {\it second formulation}:  He shows that {\bf the big exponential map not only interpolates the Bloch-Kato exponential maps themselves, but simultaneously  the dual Bloch-Kato exponential maps}, see also \cite[\S 4.2.1, especially (4.2.1)]{perrin99}. The equivalence is explained in \cite[\S 4.2]{perrin99}. We shall only state these kind of formulae later in the Lubin-Tate setting, see Theorem \ref{thm:BF} below: those formulae specialize to those on p. 123 in    Berger's article: \cite[II.16]{berger-exp} uses  the same approach for a very elegant proof based on such interpolation formulae in Thm.\ II.10 of (loc.\ cit.) calculated via the theory of  $(\varphi,\Gamma)$- and Wach-modules again.


Noting that $\Omega_{V,h} $ possesses an (almost) inverse $\mathcal{L}_{V,h}$, called {\em regulator map,}  we obtain a {\it third formulation} as follows.

\begin{thm}[Third formulation: { Adjunction of Big Exponential-  and regulator map}]\label{thm:PRIII}
The regulator map is adjoint to the big exponential map:
\begin{equation*}
\xymatrix{
   {D(\Gamma,\qp)\otimes_{\Lambda} H_{Iw}^1(V^*(1)) } \ar@{}[r]|{\times} &  {D(\Gamma,\qp)\otimes_{\Lambda} H_{Iw}^1(V)} {\ar[d]_{\mathcal{L}_{V,h}}}
   \ar[r]^-{\{,\}_{Iw}} & D(\Gamma,\qp) \ar@{=}[d] \\
      D(\Gamma,\qp)\otimes_{\qp} D_{cris}(V^*(1)) \ar[u]_{\Omega_{V^*(1),1-h}}  \ar@{}[r]|{\times} &  D(\Gamma,\qp)\otimes_\qp D_{cris}(V)\ar[r]^-{[,]_{D_{cris}} } &
      D(\Gamma,\qp).}
\end{equation*}
\end{thm}

\subsection{Iwasawa cohomology via $(\varphi,\Gamma)$-modules}\label{sec:LT}
Fontaine \cite{fontaine1990} observed that Iwasawa cohomology can be expressed in terms of $(\varphi,\Gamma)$-modules.

\parbox[t]{9cm}{  We recall his equivalence of categories directly in the Lubin-Tate setting, as introduced before Theorem \ref{thm:wiles}, see \cite{SV15} for more details:\\ $L_n=L(LT[\pi_L^n])$  is the extension of $L$  generated by the $\pi_L^n$-torsion points of $LT$, \\    $L_\infty := \bigcup_n L_n$,\\
$\Gamma_L := \Gal(L_\infty/L)$ and $\Gamma_n := \Gal( L_\infty/L_n)$;\\
 the Lubin-Tate character $\chi_{LT}$ induces  isomorphism $\Gamma_L \xrightarrow{\cong} o_L^\times$ and $Lie(\Gamma_L) \xrightarrow{\cong} L$, and $\nabla\in Lie(\Gamma_L)$ is the  preimage of $1.$ }\hspace*{0cm}
 {\parbox[t]{3cm}{$
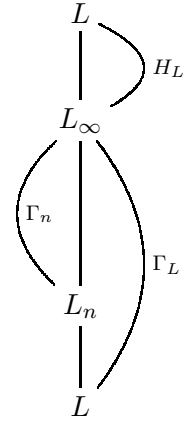
\xymatrix{
 & &\overline{L}\ar@{-}[d]\\
  & &L_\infty\ar@/_2pc/@{-}[u]_{H_L} \\
    & &\\ & & L_{n} \ar@/^2pc/@{-}[uu]_{\Gamma_n}\ar@{-}[d]\ar@{-}[uu]\\
  & &{L }\ar@/_2pc/@{-}[uuu]_{\Gamma_L} &
}$}}

 To this aim  we write
$ \mathbf{A}_L:=\widehat{o_L[[Z]][\frac{1}{Z}]}$ for the $p$-adic completion of $o_L[[Z]][\frac{1}{Z}] $. This ring comes  with commuting actions by $\varphi=[\pi_L]^*$ and $\Gamma_L$ via $[\chi_{LT}(\gamma)]^*.$ Moreover, we have the  (almost) left inverse operator $\psi_L$ satisfying $\psi_L\circ\varphi_L=\frac{q}{\pi_L},$ with $q:=\# (o_L/\pi_Lo_L)$.

Then,  $\Phi\Gamma(\mathbf{A}_L)$, the category of $(\varphi,\Gamma)$-modules, consists of finitely generated $ \mathbf{A}_L$-modules $M$ with commuting (semi-linear) actions by some endomorphism $\varphi_M$ and by $\Gamma_L$. $M$ is called \'{e}tale, if the linearized map
\begin{align*}
     \varphi_M^{lin} : \mathbf{A}_L \otimes_{\mathbf{A}_L,\varphi_L} M & \xrightarrow{\; \cong \;} M \\
                                           f \otimes m & \longmapsto f \varphi_M (m)
\end{align*} is an isomorphism. We write $\Phi\Gamma^{et}(\mathbf{A}_L)$ for the full subcategory of \'{e}tale objects, to which e.g.\ the module $\Omega^1:=\mathbf{A}_LdZ$ of differential forms belongs.   The reader is invited to  recall from the table just before Theorem \ref{thm:wiles} how the Lubin-Tate setting specializes to the cyclotomic one.

Let $\Rep_{o_L}(G_L)$ denote the abelian category of finitely generated $o_L$-modules equipped with a continuous linear $G_L$-action. The following result is established in \cite{KR} Thm.\ 1.6., where also the ring $\mathbf{A}$ is introduced.

\begin{thm}\label{KR-equiv}
The functors
\begin{equation*}
    V \longmapsto D_{LT}(V) := (\mathbf{A} \otimes_{o_L} V)^{\ker(\chi_{LT})} \qquad\text{and}\qquad M \longmapsto (\mathbf{A} \otimes_{\mathbf{A}_L} M)^{\phi_q \otimes \varphi_M =1}
\end{equation*}
are exact quasi-inverse equivalences of categories between $\Rep_{o_L}(G_L)$ and $\Phi\Gamma(\mathbf{A}_L)$.
\end{thm}

Herr observed how to express Galois cohomology in terms of $(\varphi,\Gamma)$-modules, see \cite{KV} for a generalisation in the Lubin-Tate setting. The following results, generalising Fontaine's observation from the cyclotomic to the general Lubin-Tate setting, expresses Iwasawa cohomology, where the twist by $\tau:=\chi_{LT}^{-1} \chi_{cyc}$ is a new phenomenon in the general Lubin-Tate stting (disappearing obviously in the cyclotomic case).

\begin{thm}[{\cite[Thm.\ 5.13]{SV15}}]\label{psi-version}
Let $V$  in $\Rep_{o_L}(G_L)$.   Then, with $\psi = \psi_{D_{LT}(V(\tau^{-1}))}$, we have a short exact sequence
\begin{multline}\label{f:psi-version}
  0 \longrightarrow H^1_{Iw}(L_\infty/L,V) \longrightarrow D_{LT}(V(\tau^{-1})) \xrightarrow{\;\psi-1\;}  D_{LT}(V(\tau{-1})) \\ \longrightarrow H^2_{Iw}(L_\infty/L,V) \longrightarrow 0 \ ,
\end{multline}
which is functorial in $V$.
\end{thm}

Now the above mentioned version of the {{Schmid-Witt-Formula}} \ref{schmid-witt}     for ramified Witt-vectors (of finite length)
gives an explicit determination of elements in   Iwasawa-cohomology (in Lubin-Tate setting) as images under the  Kummer map
\begin{equation*}
  \kappa : \varprojlim_{n,m} L_n^\times/{L_n^\times}^{p^m} \xrightarrow{\;\cong\;} H^1_{Iw}(L_\infty/L,\mathbb{Z}_p(1)) \ ,
\end{equation*}
by means of  $(\varphi,\Gamma_L)$-modules. We write $T^*=o_L\eta^*$ for the dual of the Tate module $T:=\varprojlim_n \mathcal{F}[\pi_L^n]=o_L\eta$ of the $\mathcal{F}$ and we consider the map
 \begin{align*}
  \nabla : (\varprojlim_n L_n^\times) \otimes_\mathbb{Z} T^* & \longrightarrow  \mathbf{A}_L^{\psi=1} ,
       u \otimes a\eta^* \longmapsto a \frac{\partial_\mathrm{inv}(g_{u,\eta})}{g_{u,\eta}}(\iota_{LT}(\eta)) \ .
\end{align*}
   \begin{thm}[{\it cycl:} {\sc Benois}{\cite[Prop.\ 2.1.5]{benois}}, {\sc Colmez}{\cite[Thm.\ 7.4.1]{colmez-tsinghua},\cite[Prop.\ V.3.2 (iii)]{cher-col1999}} | {\it LT:} {\sc Schneider-Venjakob}{\cite[Thm.\ 6.2]{SV15}}]\label{thm:Kummermap}
The diagram
  \begin{equation*}
     \xymatrix{
         (\varprojlim_n L_n^\times) \otimes_{\mathbb{Z}} T^* \ar[dr]_-{\nabla} \ar[r]^-{\kappa \otimes T^*}_-{}    &  H^1_{Iw}(L_\infty/L,o_L(\tau)) \ar[d]^-{Exp^*}_-{\cong} \\
         & {\mathbf{A}_L^{\psi=1}}=D_{LT}(o_L)^{\psi=1}  &  }
\end{equation*}
is commutative.
\end{thm}

Colmez \cite[\S 7.4]{colmez-tsinghua} calls this already an {\it explicit reciprocity law}, probably as the determination of the Kummer map is achieved by using the cup-product pairing, at least in \cite[\S 6]{SV15}. Moreover, the latter  is calculated by the explicit Schmid-Witt formula as had been noticed already by Fontaine in \S 4.4 of his letter to Perrin Riou \cite{perrin94}. Finally, as Colmez explains in Remark following \cite[Thm.\ 7.4.1]{colmez-tsinghua}
$\mathrm{Exp}^*$ produces - up to the Amice-transform and restricting (pseudo)measures from $\zp$ to $\mathbb{Z}_p^\times$  - the Kubota-Leopoldt zeta function from the system of cyclotomic units - as we discuss in detail in section  \ref{sec:ES}.\footnote{  $\nabla(\mathbf{c}(a,b))=D\log g_{\mathbf{c}(a,b)}=D\log f(T)$ with the notation in \eqref{f:logderzeta}. The restriction of measures from $\zp$ to $\mathbb{Z}_p^\times$ corresponds to applying the operator $1-\varphi=1-\varphi\psi$ given that we start with a power series fixed under $\psi$. Therefore, Colmez' interpretation is compatible with that of subsection \ref{sec:p-adicL}.} See also \cite[Thm.\ IV.2.1, Prop.\ IV.3.1 (ii)]{cher-col1999}, \cite[Thm.\ 9 (iii)]{colmez1998} for an explicit reciprocity involving $Exp^*$ for general de Rham representations interpolating dual Bloch-Kato exponential maps for the twists $V(k)$ of $V:$
\[p^{-n}\varphi^{-n}(Exp^*(\mu))=\sum_{k\in\z}\exp^*_{K_n,V(-k)}\left( \int_{\Gamma_{K_n}}\chi_{cyc}(x)^{-k}\mu(x)\right)\]
for $\mu\in H^1_{Iw}(K,V)$ and $n$ big enough.
As remarked in \cite[\S 5, Remark (ii)]{col-coleman} the term $CW_{V,k,n}(\mu)$ corresponding to $\exp^*\left( \int_{\Gamma_{K_n}}\chi_{cyc}(x)^{-k}\mu(x)\right) $ in the sum can be defined directly from $Exp^*(\mu)$ without any reference to $\exp*$ and the maps $\mu\mapsto CW_{V,k,n}(\mu)$ are generalizations of the Coates-Wiles homomorphisms from sections \ref{sec:BK} and \ref{sec:ES}. The meaning of the above formula is that they are related to Bloch-Kato's dual exponential maps, which is an explicit reciprocity law in the sense of Corollary \ref{cor:KatoExplicit}. Compare also with formulae \eqref{f:CW1} and \eqref{f:CW2} combined with the formula in Theorem \ref{thm:Kummermap}, which also indicates the relation of $Exp^*$ with Coleman's isomorphism $Col$ in \eqref{f:col-regulator} and \eqref{f:Col}.

\subsection{Colmez' abstract reciprocity law}

In \cite{col-mira} Colmez formulates a reciprocity law purely in terms of $(\varphi,\Gamma)$-modules as an abstraction of Theorems \ref{thm:PRI},\ref{thm:PRIII} above. To this end he replaces the representation $V$ by $D=D(V)\in \Phi\Gamma(\mathbf{A}_\Qp)$,  $H^1_{Iw} (V)$ by $ D^{\Psi=1}$ (which maps via $ {1-\varphi} $ to $D^{\psi=0}= D\boxtimes\mathbb{Z}_p^\times $ in his terminology),
  $ \mathcal{L}_V$ ( respectively the inverse of $\Omega_{V^*(1)} $) by $ 1-\varphi$.
Furthermore, he defines $  \mathbf{A}_\Qp(\Gamma):=\widehat{\Lambda(\Gamma)[\frac{1}{\gamma-1}]}$ as  $p$-adic completion and constructs a further pairing $ \{,\}_{convolution}$ by generalising the convolution product of distributions to $(\varphi,\Gamma)$-modules. Then his abstract reciprocity law has the following shape.

\begin{thm}[{\bf Colmez' abstract reciprocity law}]
The canonical pairing
\[\check{D}:=Hom_{\mathbf{A}_\Qp}(D,\Omega^1) \times D \to \Omega^1\] induces a commutative diagram
\[\xymatrix{
  \check{D}^{\psi=0} \ar[d]_{\sigma_{-1}\iota_*}   & \times    & D^{\psi=0} \ar@{=}[d]_{ } \ar[r]^{\{,\}_{Iw}} & \mathbf{A}_\Qp(\Gamma)   \ar[r]^{\mathfrak{M}} & \mathbf{A}_\Qp^{\psi=0}   \ar[d]^{d} \\
 \check{D}^{\psi=0}   & \times  &  D^{\psi=0}\ar[rr]^{\{,\}_{convolution}} &   & (\Omega^1)^{\psi=0} .    }\]
 Here, $\iota_*$ denotes a certain involution, $\sigma_1\in \Gamma$ corresponds to $-1$ under the cyclotomic character and $\mathfrak{M}$ denotes the Mellin-transform.
\end{thm}
This compatibility has been used by Colmez to  study  locally algebraic vectors in order to compare the $p$-adic with the classical local Langlands correspondence (for $GL_2(\Qp)$).

This abstract version in the cyclotomic setting served as a role model for an abstract reciprocity law in the Lubin-Tate setting as we will describe now. As sofar there is no good integral theory in that case we will only  obtain a version over the Robba ring.

\subsection{Abstract Reciprocity formula in the Lubin-Tate setting}\label{sec:abstract}

We start recalling the definition of the Robba ring and $(\varphi_L,\Gamma_L)$-modules in the context of Lubin-Tate extensions over the basis $L$. Let  $L\subseteq K$ be a complete field containing the period $\Omega\in K$.  The Robba ring
 $\mathcal{R}\subseteq K[[Z,Z^{-1}]]$ with coefficients in $K$ consists of those analytic functions, which converge on some annulus  $r\leq |Z|<1$  for some $0<r<1$
 \xy 0;/r5pc/:;*\dir{*}="0",*\xycircle<15pt,15pt>{},*\xycircle<11pt,11pt>{{.}}\endxy

\noindent while
 $\mathcal{R}^+:=\mathcal{R}\cap K[[Z]]$ consists of those analytic functions, which converge on the open unit disk $\mathbb{B}:=\{|Z|<1\}. $ \xy 0;/r5pc/:;*\dir{*}="0",*\xycircle<15pt,15pt>{}\endxy\\

On $\cR,\cR^+$  the elements $\gamma\in\Gamma_L$ and $\varphi_L$ act via   $Z\mapsto [\chi_{LT}(\gamma)](Z)$  and $  [\pi_L](Z)$, respectively. Moreover, there exists an   (almost) left inverse operator $\psi$ satisfying $\psi\circ\varphi_L=\frac{q}{\pi_L}.$

We write
$Rep_L(G_L)$ for the category of finite dimensional $L$-vector spaces  with continuous, linear $G_L$-action  and   $Rep_L^{an}(G_L)$ for the full subcategory consisting of {\it analytic} representations $V$,  i.e.,  $\mathbb{C}_p\otimes_{\sigma,L}V$   is   the trivial semilinear $\mathbb{C}_p$-representation   $\mathbb{C}_p^{\mathrm{dim}_L V}$  for all $\iota\neq\sigma:L\to \mathbb{C}_p$. Furthermore,
  $\Phi\Gamma^{an}({\cR_L})$ denotes the category of $(\varphi_L,\Gamma_L)$-modules\footnote{From now on we assume them to be free of finite rank over $\cR_L.$ } over $ \cR_L$ with {\bf $L$-linear} $Lie(\Gamma_L)$-action (the latter derived action is always $\Qp$-linear!).

\begin{thm}[{\sc Kisin-Ren/Fontaine, Berger}]
There is an equivalence of categories
 \begin{align*}
 Rep_L^{an}(G_L) & \longleftrightarrow \Phi\Gamma^{an,\acute{e}t}({\cR_L}) \\
  W & \mapsto D_{rig}^\dagger(W).
\end{align*}
\end{thm}
 Note that without the superscripts $\phantom{m}^{an}\phantom{m}$ the statement is false for $L\neq\Qp$, because the theorem of overconvergence by Cherbonnier and Colmez \cite{cher-col1998} does not hold, see \cite{FX} for a counterexample.

Now let   $V$ be an $L$-linear continuous representation of $G_L$ such that $V^*(1)$
 (and hence $V(\tau^{-1})$)   is $L$-analytic and crystalline. Then
 $M:=D_{\mathrm{rig}}^\dagger(V(\tau^{-1}))$ and $ \check{M}:=\Hom_{\cR}(M,\Omega^1_{\cR}) \cong D_{\mathrm{rig}}^\dagger(V^*(1))$ belong to $\Phi\Gamma^{an}(\cR)$ over the Robba ring $\cR$. Building on   {\sc Berger}'s comparison isomorphism
\[\mathrm{comp}_{\check{M}}:\cR[\frac{1}{t_{LT}}]\otimes_{\cR} \check{M}\xrightarrow{\cong} \cR[\frac{1}{t_{LT}}]\otimes_L D_{cris,L}(V^*(1))\]
 we prove the following abstract reciprocity law.
\begin{thm}[{\sc Schneider-Venjakob}{\cite[(140)]{SV23}}]
 The following diagram
\[\tiny\xymatrix@C=0.2cm{
  \check{M}^{\psi_L=0} \ar@{=}[d]_{}  \ar@{}[r]|{\times} & M^{\psi_L=0} \ar@{=}[d]_{} \ar[r]^(0.5){ \{,\}'_{M,Iw}} &\cR(\Gamma_L)  \ar[r]^{\mathfrak{M}}_{\cong} & \cR_L^{\psi=0}   \ar[d]^{d}  \\
  \check{M}^{\psi_L=0} \ar@{<-->}[d]_{\mathrm{comp}_{\check{M}}}  \ar@{}[r]|{\times} & M^{\psi_L=0} \ar@{<-->}[d]_{\mathrm{comp}_M} \ar[rr]^(0.5){} &&  (\Omega_\cR^1)^{\psi=0} \ar@{<-->}[d]_{\mathrm{comp}_{\Omega^1}}     \\
 {\cR^{\psi_L=0}\otimes_L D_{cris,L}(V^*(1))}  \ar@{}[r]|{\times} & {\cR^{\psi_L=0}\otimes_L D_{cris,L}(V(\tau^{-1}))} \ar[rr]^-{[,]_{{D_{cris,L}(V(\tau^{-1}))}}} && {\cR^{\psi_L=0}\otimes_L D_{cris,L}(L(\chi_{LT}))} } \]
 commutes up to inverting $t_{LT}:=\log_{LT}$. Here, the dotted arrows indicate that we have only an isomorphism after inverting $t_{LT}.$
\end{thm}
Actually, this theorem is  a consequence of Serre-duality on certain rigid analytic varieties comparing additive and multiplicative residuum maps. More precisely, we show  the commutativity of
\[ \xymatrix{
 \cR_L(\Gamma_L):=  \cR_L(\mathfrak{X}^\times) \ar[dd]_{(-)(\mathrm{ev}_1 d\log_\mathfrak{X})}^{\cong} \ar[r]^-{\cdot d\log_{\mathfrak{X}^\times }}_-{\cong} & \Omega^1_{\mathfrak{X}^\times} \ar[d]^{\mathrm{res}_{\mathfrak{X}^\times}} \\
       & L  \\
   (\Omega^1_{\mathfrak{X}})^{\psi=0} \ar[r]^(0.6){\mathrm{ev}_{-1}\cdot} & \Omega^1_{\mathfrak{X} } \ar[u]_{\mathrm{res}_{\mathfrak{X}}}, }\]
 with {\sc Schneider-Teitelbaum}'s character varieties $\mathfrak{X},\mathfrak{X}^\times$ for the groups $o_L, o_L^\times$ using their Fourier theory and Lubin-Tate isomorphism $\mathfrak{X}\cong \mathbb{B}$ over $\mathbb{C}_p.$ See \cite[Thm.\ 4.5.12]{SV23} for details.

From the above abstract version one derives now a reciprocity formula in the Lubin-Tate setting analogous to Theorem \ref{thm:PRIII}. We write $D(\Gamma_L,\Cp)$  for the algebra of locally $L$-analytic distributions. Recall that already Fourquaux \cite{Fou}, who initiated the investigation of Perrin-Riou's approach for  Lubin-Tate extensions in his thesis in 2005, had achieved a generalization of Colmez' construction of the {\it Perrin-Riou logarithm}. In the following result enters instead the Big exponential map, he has constructed together with Berger.

\begin{thm}\label{thm:LTRec}
If $\mathrm{Fil}^{-1}D_{cris,L}(V^*(1))=D_{cris,L}(V^*(1))$ and
$D_{cris,L}(V^*(1))^{\varphi_L=\pi_L^{-1}}=D_{cris,L}(V^*(1))^{\varphi_L=1}=0$, then the following diagram commutes:
\begin{equation*}\small
\xymatrix@=0.7em@R=1.7em{
    {D_{\mathrm{rig}}^\dagger(V^*(1))^{\psi_L=\frac{q}{\pi_L}}}\ar@{}[r]|{\times} &  {D_{LT}(V(\tau^{-1}))^{\psi_L=1}}\ar[d]_{\mathcal{L}_V^0}
    \ar[r]^(0.6){\{,\}_{Iw}} & D(\Gamma_L,\Cp) \ar@{=}[d] \\
    D(\Gamma_L,\Cp)\otimes_L D_{cris,L}(V^*(1)) \ar[u]_{\Omega_{V^*(1),1}}   \ar@{}[r]|{\times} &  D(\Gamma_L,\Cp)\otimes_L
    D_{cris,L}(V(\tau^{-1})) \ar[r]^(0.65){[,]_{}} & D(\Gamma_L,\Cp),}
\end{equation*}
where $\Omega_{V^*(1),1} $ denotes Berger's and Fourquaux' Big exponential map \cite{berger-four} while the regulator map $ \mathcal{L}_V^0$ is defined in \cite[\S 5.1]{SV23}.
\end{thm}

We write $\mathrm{Ev}_{W,n}:\mathcal{R}^+_L \otimes_L D_{cris,L}(W)\to L_n\otimes_L D_{cris,L}(W)$ for the composite $\partial_{D_{cris,L}(W)}\circ \varphi_q^{-n}$ from the introduction of \cite{berger-four}, which actually sends $f(Z)\otimes d$ to $f(\eta_n)\otimes \varphi_L^{-n}(d).$ By abuse of notation we also use $\mathrm{Ev}_{W,0}$ for the analogous map   $ \mathcal{R}^+_K \otimes_L D_{cris,L}(W)\to K\otimes_L D_{cris,L}(W)$. For $x\in  D(\Gamma_L,K)\otimes_L D_{cris,L}(W)$ we denote by $x(\chi_{LT}^{j}) $ the image under the map $ D(\Gamma_L,K)\otimes_L D_{cris,L}(W)\to K\otimes_L D_{cris,L}(W),$ $\lambda\otimes d\mapsto \lambda(\chi_{LT}^{j})\otimes d.$
With this  notation Berger's and Fourquaux' interpolation property reads as follows:

\begin{thm}[{Berger-Fourquaux \cite[Thm.\ 3.5.3]{berger-four}}]\label{thm:BF}
 Let $W$ be $L$-analytic and $h\geq 1$ such that \linebreak $\mathrm{Fil}^{-h}D_{cris,L}(W)=D_{cris,L}(W)$. For any  $f\in \left((\mathcal{R}^+)^{\psi=0}\otimes_L D_{cris,L}(W)\right)^{\Delta=0}$ and $y\in \left(\mathcal{R}^+ \otimes_L D_{cris,L}(W)\right)^{\psi=\frac{q}{\pi_L}}$ with $f=(1-\varphi_L)y$ we have: If $h+j\geq 1$, then
\begin{align}\label{f:interjgeq1}\notag
  &h^1_{L_n,W(\chi_{LT}^j)}  (tw_{\chi_{LT}^j}(\Omega_{W,h}(f)))= \\
   & (-1)^{h+j-1}(h+j-1)!\left\{
                           \begin{array}{ll}
                             \exp_{L_n,W(\chi_{LT}^j)}\Big(q^{-n}\mathrm{Ev}_{W(\chi_{LT}^j),n}( \partial_{\mathrm{inv}}^{-j}y\otimes e_j)\Big) & \hbox{if $n\geq 1$;} \\
                   \exp_{L,W(\chi_{LT}^j)}\Big((1-q^{-1}\varphi_L^{-1})\mathrm{Ev}_{W(\chi_{LT}^j),0}   ( \partial_{\mathrm{inv}}^{-j}y\otimes e_j)\Big)         , & \hbox{if $n=0$.}
                           \end{array}
                         \right.
\end{align}
If $h+j\leq  0$, then
\begin{align}\label{f:interjleq0}\notag
\exp_{L_n,W(\chi_{LT}^j)}^* &\Big(  h^1_{L_n,W(\chi_{LT}^j)}  (tw_{\chi_{LT}^j}(\Omega_{W,h}(f)))\Big)= \\
   & \frac{1}{(-h-j)!}\left\{
                           \begin{array}{ll}
                            q^{-n}\mathrm{Ev}_{W(\chi_{LT}^j),n}(\partial_{\mathrm{inv}}^{-j}y\otimes e_j) & \hbox{if $n\geq 1$;} \\
                    (1-q^{-1}\varphi_L^{-1})\mathrm{Ev}_{W(\chi_{LT}^j),0}   (\partial_{\mathrm{inv}}^{-j}y\otimes e_j)         , & \hbox{if $n=0$.}
                           \end{array}
                         \right.
\end{align}
\end{thm}

The reciprocity law Theorem \ref{thm:LTRec} then delivers the interpolation property of the regulator map
\[\mathbf{L}_V: H^1_{Iw}(L_\infty/L,T)\to  D(\Gamma_L,\Cp)\otimes_L D_{cris,L}(V(\tau^{-1})),\]
which  generalizes \cite[Thm.~A.2.3]{LVZ} and \cite[Thm.\ B.5]{zerbes-loeffler}  from the cyclotomic  case (See \S 5.2.4 of \cite[Thm.\ 6.4]{SV23} for the missing notation).
\begin{thm}[{\cite[Thm.\ 6.4]{SV23}}]\label{thm:SV} Assume that $V^*(1)$ is $L$-analytic with $\mathrm{Fil}^{-1}D_{cris,L}(V^*(1))=D_{cris,L}(V^*(1))$ \linebreak and $D_{cris,L}(V^*(1))^{\varphi_L=\pi_L^{-1}}=D_{cris,L}(V^*(1))^{\varphi_L=1}=0$. Then it holds that for   $j\geq 0$
\begin{align*}
 \Omega^{j} \mathbf{L}_V(y)(\chi_{LT}^j)&=  j!\Big((1-\pi_L^{-1}\varphi_L^{-1})^{-1}(1-\frac{\pi_L}{q}\varphi_L)\widetilde{\exp}^*_{L,V(\chi_{LT}^{-j}),\id}(y_{\chi_{LT}^{-j}})\Big)\otimes e_{ j} \\
    &=  j!(1-\pi_L^{-1-j}\varphi_L^{-1})^{-1}(1-\frac{\pi_L^{j+1}}{q}\varphi_L)\Big(\widetilde{\exp}^*_{L,V(\chi_{LT}^{-j}),\id}(y_{\chi_{LT}^{-j}})\otimes e_{ j}\Big)
\end{align*}
and for $j\leq -1:$
\begin{align*}
 \Omega^{j} \mathbf{L}_V(y)(\chi_{LT}^j)&=  \frac{(-1)^{j}}{(-1-j)!}\Big((1-\pi_L^{-1}\varphi_L^{-1})^{-1}(1-\frac{\pi_L}{q}\varphi_L)\widetilde{\log}_{L,V(\chi_{LT}^{-j}),\id}(y_{\chi_{LT}^{-j}})\Big)\otimes e_{ j} \\
    &=  \frac{(-1)^{j}}{(-1-j)!}(1-\pi_L^{-1-j}\varphi_L^{-1})^{-1}(1-\frac{\pi_L^{j+1}}{q}\varphi_L)\Big(\widetilde{\log}_{L,V(\chi_{LT}^{-j}),\id}(y_{\chi_{LT}^{-j}})\otimes e_{ j}\Big),
\end{align*}
if the operators $1-\pi_L^{-1-j}\varphi_L^{-1}, 1-\frac{\pi_L^{j+1}}{q}\varphi_L $ or equivalently $1-\pi_L^{-1}\varphi_L^{-1}, 1-\frac{\pi_L}{q}\varphi_L$ are invertible on $D_{cris,L}(V(\tau^{-1}))$  and $D_{cris,L}(V(\tau^{-1}\chi_{LT}^j))$, respectively.
\end{thm}

We leave it to the interested reader to check that the combination (of parts of) Theorem \ref{thm:BF} and \ref{thm:SV} implies again the reciprocity law \ref{thm:LTRec} by the principle of $p$-adic interpolation.

\subsection{Nakamura's Explicit Reciprocity Law}\label{sec:Nakamura}

In the cyclotomic setting Nakamura \cite{nak} extended the definition of the (dual) Bloch-Kato exponential map to all (de Rham) $(\varphi,\Gamma)$-modules over the Robba ring, i.e., they need not be \'{e}tale. In the rank one case he established an explicit reciprocity law in the style of Theorem \ref{thm:SV}. Part of his work has been carried over to the Lubin-Tate situation in \cite{MSVW}, from where we cite the following result. We define a regulator map
\[\mathbf{L}_{\cR(\delta)}: H^1_{Iw,+}(\cR(\delta)):=(\cR^+_L(\delta))^{\Psi=1}\xrightarrow{1-\varphi } (\cR^+_L(\delta))^{\Psi=0} \xleftarrow[\cong]{\mathfrak{M}_{\delta}\circ \sigma_{-1}} D(\Gamma_L), \]
where $\mathfrak{M}_{\delta}$ denotes again a certain Mellin-transform. Moreover, we recall that the Amice transform is the  map
\[A_{-}:D(o_L,K)\to \cR_K^+,\] sending a distribution $\mu$ to
\[A_\mu(Z)=\int_{o_L}\eta(x,Z)\mu(x)\]
with $\eta(x,Z):=\exp\left(\Omega x\log_{LT}(Z)\right)\in 1+Zo_{\widehat{L_\infty}}[[Z]]$.

\begin{thm}[Explicit reciprocity formula {\cite[4B1]{nak},\cite[Prop.\ 9.13/17]{MSVW}}]\label{thm:nakamura}  Let $\delta=\delta_{lc}x^k$ be de Rham.
For $k\leq 0,$ the following diagram is commutative:
\begin{equation}
\label{f:Nakamura}
\xymatrix{
  H^1_{Iw,+}(\cR(\delta)) \ar[dd]_{x\mapsto [(0,C_{Tr}({\mf Z_n})^{-1} x)] } \ar[rrr]^{\mathbf{L}_{\cR(\delta)} } &    & & D(\Gamma_L)  \ar[dd]^{pr_{\Gamma_n}} \\
      &   &  & \\
  H^1_{\Psi,{\mf Z_n}}(\cR(\delta))\ar[rr]^-{\exp^{*,(n)}}  &  & D_{dR}^{(n)}(\cR(\delta))\cong L_n\otimes_L D_{dR}(\cR(\delta))   &K[\Gamma_L/U],\ar[l]_(0.3){\Sigma} }
\end{equation}
i.e., a  class $[A_\mu \mathbf{e}_\delta]\in H^1_{\Psi,\mf Z_n}(\cR(\delta))^{\Gamma_L},$  is mapped under $\exp^{*}$ to
\[\mathfrak{C}(\delta )\int_{ o_L^\times}\delta(x)^{-1} \mu(x)\frac{1}{t_{LT}^k}\mathbf{e}_{\delta },\] where $\mathfrak{C}(\delta ) $ defines a certain $\varepsilon$-constant defined in (loc.\ cit.).
For $k\geq 1,$ the following diagram is commutative:
\begin{equation}
\label{f:Nakamura2}
\xymatrix{
 H^1_{Iw,+}(\cR(\delta))\ar[dd]_{x\mapsto [(0,C_{Tr}({\mf Z_n})^{-1} x)] }  \ar[rrr]^{\mathbf{L}_{\cR(\delta)} } &    & & D(\Gamma_L)\ar[dd]^{pr_{\Gamma_n}} \\
      &   &  & \\
  H^1_{\Psi,{\mf Z_n}}(\cR(\delta)) &   & D_{dR}^{(n)}(\cR(\delta))\cong L_n\otimes D_{dR} (\cR(\delta)) \ar[ll]_-{\exp^{(n)}} &K[\Gamma_L/U],\ar[l]_(0.3){\Sigma'} }
\end{equation}
i.e., a  class $[A_\mu \mathbf{e}_\delta]\in   H^1_{\Psi,\mf Z_n}(\cR(\delta))^{\Gamma_L},$  is mapped under $\exp^{-1}_{\cR(\delta)}$ to
\[ \mathfrak{C}'(\delta )\int_{ o_L^\times}\delta(x)^{-1} \mu(x)\frac{1}{t_{LT}^k}\mathbf{e}_{\delta },\]
where $\mathfrak{C}'(\delta ) $ defines a certain $\varepsilon$-constant defined in (loc.\ cit.).
\end{thm}

\section{Regulator maps and Euler systems}\label{sec:ES}

{\sc Coates-Wiles \cite[Introduction]{coates-wil77}:}
\begin{quote}
   The above result on cyclotomic fields [the link between the divisibility of $(2\pi i)^{-k}\zeta(k)$ by $p$ for even integers $k$ with $1<k<p-1$ and the non-triviality of the $\chi_{cyc}^k$-th eigenspace of the $G(\Qp(\mu_p)/\Qp)$-modules $U^1/C$ of $1$-units modulo the closure of cyclotomic units] was probably known, in essence, to
Kummer. However, we have been influenced by the important, and somewhat
neglected, paper ...[\cite{iw64} ]  of Iwasawa, which establishes a deeper result. Accordingly,
the present paper has been written in the spirit of Iwasawa's work. In particular,
we study explicit reciprocity laws in certain fields of division points on $E$, and
introduce analogues of the mysterious maps $\psi_n$, used by Iwasawa in the cyclotomic
case. However, since this paper was first written, a number of people have pointed
out to us that these explicit reciprocity laws can be avoided in the proof of Theorem 1 [the main result concerning BSD in that article]. While this is true ..., we have retained a discussion of them for two reasons.
Firstly, we believe that the explicit reciprocity laws form the correct conceptual
framework for discussing Kummer's $p$-adic logarithmic derivatives (which are
essential for the proof of Theorem 1). Secondly, we suspect that these laws will
play a vital role in future attempts to prove the finiteness of the Tate-Safarevic
group of $E$ when $L(E/F,1)\neq 0.$
\end{quote}

{\sc Perrin-Riou (1999)\cite[Introduction]{perrin99}:}
\begin{quote}
  Le d\'{e}veloppement de ces lois [explicite de r\'{e}ciprocit\'{e}] s'est fait en parall\`{e}le et en liaison avec le
d\'{e}veloppement de la th\'{e}orie d'Iwasawa locale; dans le cas classique, il s'agit de
l'\'{e}tude du comportement des unit\'{e}s locales sur la $\mathbb{Z}_p^\times$-extension cyclotomique $K_\infty$ \`{a} l'aide de l'application exponentielle (Iwasawa, Coates-Wiles, Coleman).
\end{quote}

\subsection{Cyclotomic untis}
\subsubsection{Special $L$-values}

The development of explicit reciprocity laws (especially from Kummer's to Iwasawa's version) happened simultaneously with the discovery of the description of special $\zeta$-values through special units. Already in Kummer's version \ref{thm:kummer} it was crucial to form "(logarithmic) derivatives" from such units, which means to express such units in terms of power series, which allow a differential calculus, and doing the manipulation that way. Inspired by Iwasawa this was immanent in the early work of Coates-Wiles already and then systematized by Coleman's invention of the power series attached to certain local units and named after him, see \ref{thm:Coleman}.

In the most easiest case to begin with, viz of the Riemann zeta-function, we have the well-known formulae of the Bernoulli numbers $\mathcal{B}_n$, $n\geq 0$,
\begin{align}\label{f:bernoullidef}
  \frac{1}{e^t-1} & =\sum_{n=0}^{\infty}\frac{\mathcal{B}_n}{n!}t^{n-1}
\end{align}
and its relation to the $\zeta$-values
 \begin{align}\label{f:bernoulli}
  \zeta(1-k) & =-\frac{\mathcal{B}_k}{k}\;\;\;\mbox{ $(k=2,4,6,\ldots).$}
\end{align}
Setting \[F(z)=\frac{e^{-\frac{a}{2}z}- e^{\frac{a}{2}z}}{e^{-\frac{b}{2}z}- e^{\frac{b}{2}z}}\]
we obtain for its logarithmic derivative
\begin{align}\label{f:logder}
  g(z) & :=\frac{d}{dz}\log F(z) \\
    & = \frac{b}{2}\bigg(\frac{1}{e^{-bz}-1}  -\frac{1}{e^{bz}-1}\bigg)-\frac{a}{2}\bigg(\frac{1}{e^{-az}-1}  -\frac{1}{e^{az}-1}  \bigg) \\
    & =\sum_{k=2\;(\mathrm{even})}^{\infty} \frac{\mathcal{B}_k z^{k-1}}{k!}(a^k-b^k)
\end{align}
where we used \eqref{f:bernoullidef} for the last equation.
Using \eqref{f:bernoulli} we thus obtain
\begin{align}\label{f:logderzeta}
  \bigg(\frac{d}{dz}\bigg)^{k-1}g(z)_{\mid z=0} & =(b^k-a^k)\zeta(1-k) \;\;\;\mbox{ $(k=2,4,6,\ldots)$}.
\end{align}
Using the transformation
\[F(z)=f(T)=\frac{(1+T)^{-\frac{a}{2}}- (1+T)^{\frac{a}{2}}}{(1+T)^{-\frac{b}{2}}- (1+T)^{\frac{b}{2}}}\;\;\;\mbox{ with } T=e^z-1\]
the differentiation $\frac{d}{dz}$ corresponds to $D:=(T+1)\frac{d}{dT}$ and we may rewrite \eqref{f:logderzeta} as
\begin{align}\label{f:trafo}
  \bigg(D^{k-1}D\log f(T)\bigg)_{\mid T=0} &=(b^k-a^k)\zeta(1-k) \;\;\;\mbox{ $(k=2,4,6,\ldots)$}.
\end{align}
For $\pi_n=\xi_n-1$ where $\xi=(\xi_n)_n\in \zp(1)$ a generator, i.e., a norm-compatible system of  $p^n$th roots of unity, the elements $f(\pi_n)=:c_n(a,b)$ are so called {\it cyclotomic units} whose importance with respect to special $\zeta$-values were discovered already by  Kummer:  the values of the
Riemann zeta function at the odd negative integers arise as the higher logarithmic
derivatives (seeing the uniformizing element as the varialbe) of them. They form a norm-compatible system $\mathbf{c}(a,b)=(c_n(a,b))$ in the inverse limit $\mathbb{U}(K_\infty):= \varprojlim_n o_{K_n}^\times$ of local units with respect to the norm maps, where $K_n:=\Qp(\xi_n)$ runs through the $p$-cyclotomic tower. In fact, $f(T)=g_{\mathbf{c}(a,b)}(T)$ is the Coleman power series attached to $\mathbf{c}(a,b) $ with regard to $\xi.$

This observation led Coates and Wiles to the general definition of the higher logarithmic derivative
homomorphism for each $k\geq 1$
\begin{align*}
  \delta_k: & \mathbb{U}(K_\infty)\to \zp \\
  \mathbf{u} & \mapsto  \bigg(D^{k-1}D\log g_{\mathbf{u}}(T)\bigg)_{\mid T=0}.
\end{align*}
In this terminology \eqref{f:trafo} now becomes
\begin{align}\label{f:CWder}
 \delta_k(\mathbf{c}(a,b)) &=(b^k-a^k)\zeta(1-k) \;\;\;\mbox{ $(k=2,4,6,\ldots)$}
\end{align}
relating {\it cyclotomic units} to $\zeta$-values.

\subsubsection{$p$-adic $L$-functions}\label{sec:p-adicL}
We will now extend this   relation  to the Kubota-Leopoldt $p$-adic $\zeta$-function. More precisely, for $k\geq 1, $ there is a commutative diagram (\cite[3.5.2]{co-su-buch-MC}
\[\xymatrix{
   \mathbb{U}(K_\infty) \ar[d]_{\delta_k} \ar[r]^{Col} & {\zp}[[\Gamma]] \ar[d]^{\chi_{cyc}^k} \\
  {\zp}\ar[r]^{1-p^{k-1}} &{\zp}   }\]
with $Col$ being the composite of
\[\mathbb{U}(K_\infty)\to \zp[[T]]^{\psi=0}, \mathbf{u}   \mapsto \frac{1}{p}\log\left(\frac{g_{\mathbf{u}}^p}{\varphi(g_{\mathbf{u}})}\right),\]
with the inverse of the Mellin transform
\[\mathfrak{M}:{\zp}[[\Gamma]]\xrightarrow{\cong} \zp[[T]]^{\psi=0}, \lambda\mapsto \lambda\cdot(1+T)\;\;\mbox{ with } \gamma\cdot (1+T)=(1+T)^{\chi_{cyc}(\gamma)},\]
and where the integration map $\chi_{cyc}^k:{\zp}[[\Gamma]]\to \zp$ sends $\gamma$ to $\chi_{cyc}^k(\gamma).$
We shall also write $\lambda(\chi):=\chi(\lambda)$. We obtain that $Col(\mathbf{c}(a,b) )$ is the measure satisfying
\[Col(\mathbf{c}(a,b) )(\chi_{cyc}^k)=(b^k-a^k)(1-p^{k-1})\zeta(1-k)\;\;\;\mbox{ $(k=2,4,6,\ldots)$},\] i.e., the Kubota-Leopoldt $p$-adic $\zeta$-function is nothing else than
\begin{equation}\label{f:p-adicL}
  \frac{Col(\mathbf{c}(a,b) ) }{\sigma_b-\sigma_a}
\end{equation}
where $\chi_{cyc}(\sigma_a)=a$.

$Col$ is the prototype of a regulator map! Indeed, the {\it Perrin-Riou regulator map} for crystalline representations $V$ defined by Lei-Loeffler-Zerbes in \cite{lei-loeffler-zerbes-ColemanMap,lei-loeffler-zerbes-Wach}
\[H^1_{Iw}(V) \xrightarrow{\mathcal{L}_{V,1}} D(\Gamma,\Qp)\otimes_\Qp D_{cris}(V)\]
fits into the following commutative diagram
\begin{equation}\label{f:col-regulator}
   \xymatrix{
    \mathbb{U}(K_\infty) \ar[d]_{Kummer} \ar[r]^{Col} & {\zp[[\Gamma]]} \ar[d]^{incl\otimes e_1 } \\
   { H^1_{Iw}(\Qp(1)) } \ar[r]^-{\ell_0^{-1}\mathcal{L}_{\Qp(1),1}} &  {D(\Gamma,\Qp)\otimes_{\Qp} D_{cris}(\Qp(1)) }}
\end{equation}
which is a immediate consequence of Theorem \ref{thm:Kummermap}  and the definitions.\footnote{Compare with the diagram in \cite[Appendix C]{LVZ} and with \cite[(4) and p.\ 2402]{ven-eps}} Here $\ell_0$ denotes the element in the locally $\Qp$-analytic distribution algebra $D(\Gamma,\Qp)$ of $\Gamma$ induced by the element $\nabla\in Lie(\Gamma)$ according to \cite[\S 2.3]{Schneider-Teitelbaum}, see also \cite[Def.\ 3.3.3]{LVZ}. Here $e_1=t^{-1}\otimes t_1$ denotes the natural $\Qp$-basis of $ D_{cris}(\Qp(1))$ with $t\in B_{cris}$ Fontaine's $p$-adic period and $t_1$ a $\zp$-basis of $\zp(1)$.

\subsection{Elliptic units I}
\subsubsection{Special $L$-values}
Now we consider the setting of  \cite{coates-wil78} and let $K$ be an imaginary quadratic field with class number $1$, and $\cO$
the ring of integers of $K$. Let $E$ be an elliptic curve defined over $K,$ with complex
multiplication by $\cO$, and let $\psi$ be the Gr\"{o}ssencharacter of $E$ over $K$.   Choose $p>3$ to be a
rational prime, at which $E$ has good reduction at all primes of $K$ above $p$, and which splits in $K,$
say into $(p) = \frak{p}\overline{\frak{p}}.$ We put $\pi = \psi(\frak{p}),$ so that $\pi$ is a
generator of $\frak{p}.$ Now, as $E$ has complex multiplication by $\cO$, we can also view $\pi$
as an endomorphism of $E.$ For each $n\geq 1,$ , let $E_{\pi^n}$ be the kernel of the endomorphism $\pi^n$  of $E$. We set  $F_n = K(E_{\pi^n})$ and $F_\infty=\bigcup_n F_n.$   Let $U_n$ be the local
principal units, i.e., which are congruent $1 \mod \frak{p}_n,$ of the completion of $F_n$ at the uniqe prime $\frak{p}_n$ above $\frak{p}$.

  Let $\cL$ be the period lattice of the Weierstrass $\wp$-function
associated with a Weierstrass model for $E$. Since $K$ has class number $1,$ there exists $\Omega\in \cL$, such
that $\cL = \Omega\cO$. For each integer $k\geq 1,$ let $L(\overline{\psi}^k,s)$ be the complex Hecke $L$-function
of $\overline{\psi}^k$. It has been shown by Hurwitz, Birch and Swinnerton-Dyer,
and Damerell that the numbers $\Omega^{-k}L(\overline{\psi}^k,k), \;k= 1,2,...\; ,$ belong to $K.$ We can
therefore view these numbers as lying not only in the complex field, but also in
the completion $K_\frak{p}$ of K at the non-archimedean prime $\frak{p}.$

  Let $\kappa:\Gamma \to \mathbb{Z}_p^\times  $ be the canonical character giving the action of $\Gamma=G(F_\infty/K)$ on the Tate module $T_\pi E$.    Write  $\frak{f}$ for the conductor of $\psi$ and fix a generator $f\in K$ of it.  For each integer $k\geq 1,$ write
\[\mu_k=12(-1)^{k-1}(k-1)!\left( \frac{\Omega}{f} \right)^{-k}.\]
 We consider pairs $\frak{s}=(\{\frak{a}_j\mid j\in J\},\{n_j\mid j\in J\})$ of sets consisting of integral ideals $\frak{a}_j$ prime to $6\frak{fp}$ and integers $n_j$, respectively, over an arbitrary finite index set $J$ satisfying $\sum_J n_j(N\frak{a}_j-1)=0$ for the absolute norms $N\frak{a}_j.$ Attached to such $\frak{s}$ comes
 \[\Theta(z,\frak{s}):=\prod_J \Theta(z,\frak{a}_j)^{n_j},\]
where $\Theta(z,\frak{a})$ is an elliptic function for the lattice $\cL$ defined in (loc.\ cit.), and for each integer $k\geq 0$
 \[h_k(\frak{s}):=\sum_J n_j(N\frak{a}_j-\psi^k(\frak{a}_J)).\] Moreover,
 \[\Lambda(z,\frak{s}):=\prod_J \Lambda(z,\frak{a}_j)^{n_j}\]
 with
 \[\Lambda(z,\frak{a}):=\prod_{\frak{b}\in B} \Theta(z+\psi(\frak{b})\frac{\Omega}{f},\frak{a}),\]
 where $B$ denotes a set of representatives, which are integral and prime to $\frak{f},$ of the ray class group of $K$ modulo $\frak{f}.$

The Bernoulli numbers $\mathcal{B}_n$ from our  previous example are now replaced by the Eisenstein numbers $\mathcal{C}_k(\frak{a}),$ i.e., the coefficients of the Eisenstein series
attached to the logarithmic derivative
\begin{align}\label{f:eisensteindef}
 g(z):= \frac{d}{dz}\log \Lambda(z,\frak{s}) & =\sum_{n=1}^{\infty}\mathcal{C}_k(\frak{s})z^{n-1}
\end{align}
 and its relation to the $L$-values
 \begin{align}\label{f:eisenstein}
  \mathcal{C}_k(\frak{s}) & = 12(-1)^{k-1}\left(\frac{\Omega}{f}\right)^{-k}h_k(\frak{s})L(\overline{\psi}^k,k)\;\;\;\mbox{ $(k\geq 1).$}
\end{align}
We thus obtain
\begin{align}\label{f:logderL}
  \bigg(\frac{d}{dz}\bigg)^{k-1}g(z)_{\mid z=0} & =(k-1)!\mathcal{C}_k(\frak{s})=\mu_kh_k(\frak{s}) L(\overline{\psi}^k,k)\;\;\;\mbox{ $(k\geq 1).$}
\end{align}
Using the transformation
\[\Lambda(z,\frak{s})=R(T,\frak{s}) \;\;\;\mbox{ with } T=\exp_{\hat{E}}(z)\mbox{ and } z=\log_{\hat{E}}(T)\]
the exponential and logarithm, respectively, map of the formal Lubin-Tate group $\hat{E}$, the differentiation $\frac{d}{dz}$ corresponds to $D:=\frac{1}{\log'_{\hat{E}}}\frac{d}{dT}$ and we may rewrite \eqref{f:logderL} as
\begin{align}\label{f:trafoL}
  \bigg(D^{k-1}D\log R(T,\frak{s})\bigg)_{\mid T=0} &=\mu_kh_k(\frak{s}) L(\overline{\psi}^k,k)\;\;\;\mbox{ $(k\geq 1).$}
\end{align}
For   $u=(u_n)_n\in T_\pi \hat{E}$ a basis as $\zp$-module and $\tau_n\in \mathbb{C}/\cL$ with \[\exp_{\hat{E}}(\tau_n)=-2\frac{\wp(\tau_n)}{\wp'(\tau_n)}=-2\frac{x(\tau_n)}{y(\tau_n)}=u_n,\] the elements $\Lambda(\tau_n,\frak{s})=:e_n(\frak{s})$ are so called {\it elliptic units} (due to Robert (1973)) whose importance with respect to special $L$-values is the following: the values of the
Hecke $L$-function at the positive integers arise as the higher logarithmic
derivatives (seeing the uniformizing element as the variable) of them. They form a norm-compatible system $\mathbf{e}(\frak{s})=(e_n(\frak{s}))$ of global units in $F_n$ inducing also norm-compatible elements in the inverse limit $\mathbb{U}(L_\infty):= \varprojlim_n U_n$ of principal local units of $L_n$ with respect to the norm maps, where $ L_n:=(F_n)_{\frak{p}_n}=\Qp(\hat{E}[\pi_n])$ runs through the Lubin-Tate tower attached to $\hat{E}$. In fact, $R(T,\frak{s})=g_{\mathbf{e}(\frak{s})}(T)$ is the Coleman power series with coefficients in $\zp$ attached to $\mathbf{e}(\frak{s})$ with regard to $u:$
\[R(u_n,\frak{s})=\Lambda(\tau_n,\frak{s})\]
by \cite[Theroem 5 and (16)]{coates-wil78}.

Again Coates and Wiles considered a higher logarithmic derivative
homomorphism for each $k\geq 1$
\begin{align*}
  \delta_k: & \mathbb{U}(L_\infty)\to \zp \\
  \mathbf{u} & \mapsto  \bigg(D^{k-1}D\log g_{\mathbf{u}}(T)\bigg)_{\mid T=0}.
\end{align*}
In this terminology \eqref{f:trafoL} now becomes
\begin{align}\label{f:CWderelliptic}
 \delta_k(\mathbf{e}(\frak{s})) &=\mu_kh_k(\frak{s}) L(\overline{\psi}^k,k)\;\;\;\mbox{ $(k\geq 1)$}
\end{align}
relating {\it elliptic units} to Hecke $L$-values, see \cite[(21)]{coates-wil78}.

 This relation between elliptic units and special $L$-values was important in the works of Coates and
Wiles (\cite{coates-wil77}) on Birch and Swinnerton-Dyer conjectures, see also \cite{deSh}.

\subsubsection{$p$-adic $L$-functions}
We will now extend this   relation  to the   $p$-adic $L$-function. We fix an isomorphism
\[\eta:\hat{\mathbb{G}}_m\to \hat{E};\;\;\; T=\eta(S)=\Omega_p S+ \cdots\in\widehat{{\zp}^{nr}}[[S]].\]   Then, for $k\geq 1, $ there is a commutative diagram (\cite[I\S 3.5 (11)]{deSh}
\[\xymatrix{
   \mathbb{U}(L_\infty) \ar[d]_{\delta_k} \ar[r]^{Col} & {\widehat{{\zp}^{nr}}[[\Gamma]] } \ar[d]^{\kappa^k} \\
  {\zp}\ar[r]^{\Omega_p^k(1-\frac{\pi^{k}}{p})} &{\widehat{{\zp}^{nr}}}   }\]
with $Col$ being now the composite of
\[\mathbb{U}(L_\infty)\to \zp[[T]]^{\psi_{\hat{E}}=0}, \mathbf{u}   \mapsto \frac{1}{p}\log\left(\frac{g_{\mathbf{u}}^p}{\varphi_{\hat{E}}(g_{\mathbf{u}})}\right),\]
\[\zp[[T]]^{\psi_{\hat{E}}=0}\to \widehat{{\zp}^{nr}}[[S]]^{\psi_{\hat{\mathbb{G}}_m}=0}, T\mapsto \eta(S),\]
and the inverse of the Mellin transform
\[\mathfrak{M}:{\widehat{{\zp}^{nr}}}[[\Gamma]]\xrightarrow{\cong} \widehat{{\zp}^{nr}}[[S]]^{\psi_{\hat{\mathbb{G}}_m}=0}, \lambda\mapsto \lambda\cdot(1+T)\;\;\mbox{ with }  \gamma\cdot (1+T)=(1+T)^{\kappa(\gamma)},\]
and where the integration map $\kappa^k:{\widehat{{\zp}^{nr}}}[[\Gamma]]\to \widehat{{\zp}^{nr}}$ sends $\gamma$ to $\kappa^k(\gamma).$
We shall also write $\lambda(\kappa^k):=\kappa^k(\lambda)$.

 We obtain that $Col(\mathbf{e}(\frak{s}) )$ is the measure satisfying
\[Col(\mathbf{e}(\frak{s}) )(\kappa^k)=\Omega_p^k(1-\frac{\psi(\frak{p})^{k}}{N\frak{p}})\mu_kh_k(\frak{s}) L(\overline{\psi}^k,k)\;\;\;\mbox{ $(k\geq 1)$},\] i.e., the   $p$-adic $L$-function is nothing else than
\begin{equation}\label{f:p-adicL-elliptic}
 \frac{Col(\mathbf{e}(\frak{s}) ) }{h(\frak{s})}
\end{equation}
where $h(\frak{s})$ is the measure which interpolates the factor $ h_k(\frak{s}) $ (its existence follows from the discussion before \cite[Thm.\ 18]{coates-wil78}).

\subsection{Elliptic Units II}\label{sec:reg-ex}

Now we consider again the general Lubin-Tate situation as in section \ref{sec:LT}.

Setting as before $\mathbb{U}:=\mathbb{U}(L_\infty):= \varprojlim_n  o_{L_n}^\times $ with   transition maps given by the norm we are looking for a map
\[\mathcal{L}: \mathbb{U}\otimes_{\mathbb{Z}}\TLT^* \to D(\Gamma_L,\Cp)\otimes_L D_{cris,L}(L(\tau) )\]

such that
\begin{equation}\label{f:CWinterpolation}
  \frac{ \Omega^r}{r!} \frac{1-\pi_L^{-r}}{1-\frac{\pi_L^r}{q}} \mathcal{L}(u\otimes a \eta^*)(\chi_{LT}^{ r})\otimes (t_{LT}^{r-1} \otimes \eta^{\otimes -r+1} )=CW(u\otimes a \eta^{\otimes -r})
\end{equation}
for all $r\geq 1,u\in \mathbb{U},a\in o_L,$ where $CW$ denotes the diagonal map in Corollary \ref{cor:KatoExplicit}.
We set $\mathcal{L}=\frak{L}\otimes \mathbf{d}_1$ with $\frak{L}$ given as follows
\[\frak{L}: \mathbb{U}\otimes \TLT^*\xrightarrow{\nabla} o_L[[\omega_{LT}]]^{\psi_L={1}}\xrightarrow{(1-\frac{\pi_L}{q}\varphi)}{\mathcal{O}_{\Cp}(\mathbf{B})}^{\psi_L=0}\xrightarrow{\log_{LT}\cdot}{\mathcal{O}_{\Cp}(\mathbf{B})}^{\psi_L=0}\xrightarrow{\mathfrak{M}^{-1}}D(\Gamma_L,\Cp),\]
where the map $\nabla$ has been defined before Theorem \ref{thm:Kummermap} as
the homomorphism
\begin{align*}
  \nabla : \mathbb{U} \otimes_\mathbb{Z} T^* & \longrightarrow  o_L[[\omega_{LT}]]^{\psi=1} \\
       u \otimes a\eta^* & \longmapsto a \frac{\partial_\mathrm{inv}(g_{u,\eta})}{g_{u,\eta}}(\omega_{LT}) \ .
\end{align*}
 Note that due to the multiplication by $\log_{LT}$ the maps $\mathcal{L},\; \frak{L}$ are not $\Gamma_L$-equivariant.
In \cite[\S 5.1.1]{SV23} it is shown that
%
indeed $\mathcal{L}$ satisfies \eqref{f:CWinterpolation} and that
%
%
%
%
\begin{equation}
 \label{f:claim}\mathbb{U}\otimes_{\mathbb{Z}} \TLT^*  \xrightarrow{-\kappa \otimes \TLT^*}  H^1_{Iw}(L_\infty/L,o_L(\tau)) \xrightarrow{\mathcal{L}_{L(\tau\chi_{LT})}\otimes o_L(\chi_{LT}^{-1}) \otimes t_{LT}}D(\Gamma_L,\Cp)\otimes_L D_{cris,L}(L(\tau))
\end{equation}
  coincides with
 \[\mathcal{L}: \mathbb{U}\otimes_{\mathbb{Z}}\TLT^* \to D(\Gamma_L,\Cp)\otimes_L D_{cris,L}(L(\tau) ).\]
We set $\mathfrak{l}_i := t_{LT} \partial_\mathrm{inv} - i$, $\partial_\mathrm{inv}= \frac{d}{dt_{LT}}$ and
by abuse of notation we   also write $\mathfrak{l}_i=\nabla_{\mathrm{Lie}}-i$ for the corresponding element in  $D(\Gamma_L,K).$

If one defines  $\mathcal{L}_{L(\tau)}$ \footnote{Since the representation $L(\tau)$ does not satisfy the conditions for the definition of the regulator map at the beginning of chapter 5.1 in  (loc.\ cit.) while $ L(\tau\chi_{LT}) $ does.} as a  twist of $\mathcal{L}_{L(\tau\chi_{LT})}$  by requiring the commutativity of the following diagram:
 \begin{equation*}
   \xymatrix{
     H^1_{Iw}(L_\infty/L,o_L(\tau))  \ar[d]_{\cong } \ar[rr]^-{ \mathcal{L}_{L(\tau)}} && D(\Gamma_L,\mathbb{C}_p)\otimes_L D_{cris,L} (L(\tau)) \ar[d]^{\frac{\nabla_{\mathrm{Lie}}Tw_{\chi^{-1}}}{\Omega} \otimes t_{LT}^{-1} }\\
      H^1_{Iw}(L_\infty/L,o_L(\tau\chi_{LT}))\otimes_{o_L} o_L(\chi_{LT}^{-1}) \ar[rr]^-{\mathcal{L}_{L(\tau\chi_{LT})}\otimes o_L(\chi_{LT}^{-1}) } && D(\Gamma_L,\mathbb{C}_p)\otimes_L D_{cris,L}(L(\tau\chi_{LT}))\otimes_L  L(\chi_{LT}^{-1}), }
 \end{equation*}
 then
\[\mathcal{L}: \mathbb{U}\otimes_{\mathbb{Z}}\TLT^* \to D(\Gamma_L,\Cp)\otimes_L D_{cris,L}(L(\tau) )\] also coincides with
\begin{equation}
\label{f:regulatortau}
\mathbb{U}\otimes_{\mathbb{Z}} \TLT^*  \xrightarrow{-\kappa \otimes \TLT^*}  H^1_{Iw}(L_\infty/L,o_L(\tau)) \xrightarrow{(\frac{1}{\Omega } \nabla_{\mathrm{Lie}}Tw_{\chi^{-1}}\otimes\id)\circ\mathcal{L}_{L(\tau)}}D(\Gamma_L,\Cp)\otimes_L D_{cris,L}(L(\tau)).
\end{equation}

\begin{example}
    We refer the interested reader to  \S 5 of \cite{ST-fourier} for an example of a CM-elliptic curve $E$ with supersingular reduction  at $p$ in which they attach to a norm-compatible sequence of elliptic units $e(\mathfrak{a})$ (in the notation of \cite[II 4.9]{deSh}) a   distribution $\mu(\mathfrak{a})\in D(\Gamma_L,K)$ in \cite[Prop.\ 5.2]{ST-fourier} satisfying a certain interpolation property with respect to the values of the attached (partial) Hecke-$L$-function. Without going into any detail concerning their setting and instead referring the reader to the notation in (loc.\ cit.) we just want to point out that up to twisting this distribution is the image of $\kappa(e(\mathfrak{a}))\otimes \eta^{-1}$ under the regulator map $\mathcal{L}_{L(\tau)}$:
 \[\mathcal{L}_{L(\tau)}(\kappa(e(\mathfrak{a}))\otimes \eta^{-1})=\Omega Tw_{\chi_{LT}}(\mu(\mathfrak{a}))\otimes \mathbf{d}_1.\]
 Here, $L=\mathbf{K}_p=\mathbf{F}_\wp$ (in their notation) is the unique unramified extension of $\mathbb{Q}_p$ of degree $2,$ $\pi_L=p,$ $q=p^2,$ and the Lubin-Tate formal group is $\hat{E}_\wp$ while $K=\widehat{L_\infty}$.

 Indeed, we have a commutative diagram
 \begin{equation}\label{f:Col}\xymatrix{
   \mathbb{U} \ar[d]_{-\kappa(-)\otimes \eta^{-1}} \ar[rr]^-{Col} && D(\Gamma_L,K) \ar[d]^{\Omega Tw_{\chi_{LT}}\otimes \mathbf{d}_1} \\
   H^1_{Iw}(L_\infty/L, o_L(\tau)) \ar[rr]^-{\mathcal{L}_{L(\tau)}} &&  D(\Gamma_L,K)\otimes_L D_{cris,L}(L(\tau) ),   }
 \end{equation}
 where the Coleman map $Col$ is given as the composite in the upper line of the following commutative diagram
 \begin{equation}\label{f:ColL}
   \xymatrix{
     \mathbb{U} \ar[d]_{} \ar[r]^{\log g_{-}} & \cO^{\psi_L=\frac{1}{\pi_L} } \ar[d]_{\partial_\mathrm{inv}} \ar[r]^{1-\frac{\varphi_L}{p^2} } & \cO^{\psi_L=0 } \ar[d]_{\partial_\mathrm{inv}} \ar@{=}[r]^{} & \cO^{\psi_L=0 } \ar[d]_{\mathfrak{l}_0} \ar[r]^{\mathfrak{M}^{-1}} & D(\Gamma_L,K) \ar[d]^{\nabla_\mathrm{Lie}} \\
     \mathbb{U}\otimes T^*_p \ar[r]^{\nabla} &  {\cO^{\psi_L=1 }} \ar[r]^{1-\frac{\pi_L}{q}\varphi_L } & {\cO^{\psi_L=0 } } \ar[r]^{  \log_{LT}\cdot} & {\cO^{\psi_L=0 }}  \ar[r]^{\mathfrak{M}^{-1}} & D(\Gamma_L,K), }
 \end{equation}
 in which the second line is just $\mathfrak{L}.$ Then the commutativity of \eqref{f:Col} follows by comparing \eqref{f:ColL} with \eqref{f:regulatortau}. Finally, $Col(e(\mathfrak{a}))=\mu(\mathfrak{a})(=\mathfrak{M}^{-1}(g_\mathfrak{a}(Z))\mbox{ in their notation)}$ holds by construction in (loc.\ cit.) upon noting that  on ${{\mathcal{O}}^{\psi_L=\frac{1}{\pi_L}} }$ the operator $1-\frac{\pi}{p^2}\varphi_L\circ\psi_L$, which is used implicitly to define $g_\mathfrak{a}(Z)(= (1-\frac{\pi}{p^2}\varphi_L\circ\psi_L )\log Q_\mathfrak{a}(Z))$, equals $1-\frac{\varphi_L}{p^2} .$
\end{example}

 \begin{example}\label{Manji}
Another application of $L$-analytic regulator map is about to show up in the ongoing PhD project of   Muhammad Manji \cite{manji} supervised by David Loeffler.  Just like the Perrin-Riou regulator map above, one hopes that the $L$-analytic regulator map has global applications to Iwasawa theory, in particular that it can see $p$-adic $L$-functions which lie in the $L$-analytic distribution algebra. The below example takes again $L/\mathbb{Q}_p$ unramified quadratic, where one already sees rich new behaviour and encounters a lot of obstacles. He is considering a $p$-adic Galois representation $V$ associated to an ordinary automorphic representation $\Pi$ defined over the unitary group $GU(2,1)$ with respect to the imaginary quadratic field $\mathbf{K}$, assuming $p$ is inert in $\mathbf{K}$. As in the previous example we consider  $L=\mathbf{K}_p $, $\pi_L=p,$ $q=p^2.$ 

	Consider the $2$-variable Iwasawa algebra $\Lambda=\mathcal{O}_L\kll\mathcal{O}_L^\times\krr$ and the $L$-analytic distribution algebra   $\Lambda^L_\infty = D(\mathcal{O}_L^\times,K) \cong \Gamma(\mathscr{W}^L,\mathcal{O}_{\mathscr{W}^L})$, where $\mathscr{W}^L$ is the locus of $L$-analytic characters in $\mathrm{Spf}(\Lambda)^{\mathrm{rig}}$. This is the ring Manji expects the $p$-adic $L$-function for $V$ to lie in, although it has not been constructed yet. For the lack of such a construction he takes $L_p^*(V) := \mathcal{L}_{S}(c^\Pi)$ where $S$ is a $1$-dimensional subquotient of $V$ and $c^\Pi$ is the Loeffler--Skinner--Zerbes Euler system for $\Pi$  constructed in \cite{LSZ-GSP4-EMS}. He must impose the condition that $S^*(1)$ is $L$-analytic.  The explicit reciprocity law in this case would then be to show that $L_p^*(V)$ is really ''the'' $p$-adic $L$-function $L_p(V)$ for $V$. Moreover, one  hopes that such $L_p(V)$ would be bounded, but the notion of boundedness in $\Lambda^L_\infty$ is not so well understood, see \cite{ard-ber} regarding this unsolved problem.
	\par There are some hitches in the theory, for example Iwasawa cohomology base changed to $\Lambda^L_\infty$ no longer compares with $\psi=1$ invariants of the associated $L$-analytic $(\phi,\Gamma)$-module when $L\neq \mathbb{Q}_p$, see \cite{steingart,steingart2022iwasawa}. With some adjustments Manji can define a local condition at $p$ which he uses to define a Selmer group $\widetilde{H}^i(V) \subset H^i_{Iw}(\mathbf{K}_{\Sigma_p}/\mathbf{K},V) \otimes \Lambda^L_\infty $ where $\mathbf{K}_{\Sigma_p}$ denotes the maximal abelian extension of $\mathbf{K}$ unramified away from $p$. From this he  is about to obtain a statement towards a main conjecture:
\begin{thm}[\cite{manji}] Suppose $V$ satisfies some technical conditions and that $L_p^* \neq 0$. Then, \[\mathrm{char}_{\Lambda_\infty^L} \left(  \widetilde{H}^2(V)\right) \bigg| \left( L_p^*(V) \right).\]
\end{thm}
 	The  full statement of an Iwasawa main conjecture in this context would be the following:
 \begin{conj} Under (possibly relaxed) hypothesis,
 		\[\mathrm{char}_{\Lambda_\infty^L} \left(  \widetilde{H}^2(V)\right) = \left( L_p(V) \right).\]
 \end{conj}
\end{example}	

\subsection{Kato's explicit reciprocity law}

 {\sc Kato \cite[Introduction]{kato-explicit}, \cite{kato-generalized}:}
\begin{quote}
The explicit reciprocity law is classically a mysterious relation between
Hilbert symbols and differential forms.

\noindent
The classical explicit reciprocity    [....] is concerned with the explicit description of this homomorphism [$\lambda_m $ \eqref{f:lambdam}] by using differential forms.
 \end{quote}

 In (one of the main applications of) \cite{kato-generalized} Kato replaces the tower $K_n$ of number fields by an tower of (open) modular curves $Y_n$ and he forms the inverse limit $\varprojlim_n K_2(Y_n)$ with respect to the norm maps instead of $\varprojlim_n K_1(K_n)$. Then he defines   maps $\lambda_m$, which look like the following when specialized to this application
 \begin{equation*}
   \lambda_m: \varprojlim_n K_2(Y_n)\to \varprojlim_n H^1(Y(m),\z/p^n\z(1))\to H^1(\Qp,V)\xrightarrow{\exp^*_{\Qp,V}} M_2(X_m)\otimes \Qp
 \end{equation*}
 with $V=H^1(Y(m)_{\bar{\mathbb{Q}}_p},\Qp)(1)$ and $M_2(X_m)$ being the space of weight $2$ modular forms on the compact modular curve $X_m$.  These maps generalize \eqref{f:lambdacompBK}, \eqref{f:lambdacomp}.

 The explicit reciprocity laws Theorem 4.3.1/4 (and 6.1.9) in (loc.\ cit.) express these maps $\lambda_m$ in terms of differential forms. They contain Theorem \ref{thm:wiles}, Theorem \ref{thm:BK} and Corollary \ref{cor:KatoExplicit} as special cases.

 The Euler system in this case is given by Beilinson elements in $K_2$ of modular cuves, which are analytically related to  $\lim_{s\to 0}s^{-1}L(f,s)$ for elliptic cusp forms $f$ of weight two via regulator maps while via the maps $\lambda_m$ to the values $L(f,k-1)$.

 See \cite{li} for a more detailed explanation.

 \subsection{Euler systems and localisation at $\ell= p$} \label{sec:l=p}
In this last subsection we sketch the general philosophy.
Morally, Iwasawa cohomology groups are the groups where the local parts at $p$
of Euler systems of a global $p$-adic representation $V$ live.   Perrin-Riou's regulator map $\mathcal{L}_V$ (often also called big logarithm map $Log_V$) sends them $p$-adic $L$-functions considered as measures or distributions.

In a diagram:
 \begin{align*}
   H^1_{Iw}(\Q, V )\xrightarrow{loc_p} H^1_{Iw}(\Qp, V )\xrightarrow{ Log_V} & \mbox{Measures/Distributions}\\
    \mbox{$p$-adic Euler systems} \mapsto\phantom{mmmmmm} &\mbox{$p$-adic $L$-functions .}
 \end{align*}
See \cite{col-Lpad,colmez-BSD-bourbaki} for further references on this subject and  the survey \cite{BCDDPR} for a further discussion and more examples of Euler systems.

 Recently there have been new examples by  Loeffler-Zerbes et. al. \cite{KLZ,LLZ-ES-Hilbert,LLZ-ES-imag,LPSZ-Hida,LSZ-GSP4-EMS,LSZ-GU,LSZ-syntomic,LZ-BSD,LZ-ES-local,LZ-rankin,LZ-symmetricsquare,BLLV}  and Bertolini-Darmon et. al.  \cite{DR-ES,DR-diagonal,DR-Kato,DL,BD-Kato,BDV,BBV,BDR-BFI,BDR-BFII,BSV,BSV-BK} on generalized Kato classes and progress on the BSD-conjecture.

 \subsection{Heegner points and localisation at $\ell\neq p$}\label{sec:lneqp}

There is yet another class of statements which are referred to as ''reciprocity laws'' as was pointed out to us by David Loeffler: these are the reciprocity laws in the theory of level-raising congruences for Heegner points going back to Bertolini and Darmon's paper \cite{BD05}, and more recently back in fashion again with the work of Yifeng Liu et al. These "reciprocity laws" are describing the localisation $loc_\ell$ of $p$-adic Euler systems at primes $\ell \neq p$, in contrast to the reciprocity laws for Euler systems described in subsection \ref{sec:l=p}, which are about describing the localisation at $\ell = p.$

Without going into too much details the setting in (loc.\ cit.) is as follows: $K$ denotes an imaginary quadratic field and $K_\infty$ its anti-cyclotomic $\mathbb{Z}_p$-extension, i.e., the unique $\mathbb{Z}_p$-extension of $K$ such that $c\gamma c=\gamma^{-1}$ for all $\gamma\in G_\infty=G(K_\infty/K)\cong \zp$ for complex conjugation $c.$ Then $\Lambda=\zp[[G_\infty]]$ denotes the corresponding Iwasawa algebra. Attached to certain ordinary eigenforms $f$ (on a certain finite graph $\mathcal{T}/\Gamma$ attached to a certain definite quaternion algebra $B$) they construct an element $\mathcal{L}_f\in \Lambda$ such that $L_p(f,K):=\mathcal{L}\mathcal{L}^\iota$ - with the involution $\iota$ on $\Lambda$ sending $\gamma$ to its inverse - denotes the anti-cyclotomic $p$-adic Rankin-$L$-function attached to $f$. Using Heegner points the authors construct elements $\kappa(\ell)$ in global Iwasawa cohomology groups $H^1_{Iw}(K_\infty, T_{f}/p^nT_f)$ indexed by so called $n$-admissible primes $\ell$ attached to $f$ (necessarily prime to $p$), where $T_f\subseteq V_f$ denotes a Galois stable $\zp$-lattice in the Galois representation $V_f$ attached to $f.$ They also consider semi-local Iwasawa cohomology groups together with a decomposition of $\Lambda$-modules
\[H^1_{Iw}(K_{\infty,\ell},T_f/p^nT_f)\cong H^1_{Iw, fin}(K_{\infty,\ell},T_f/p^nT_f)\oplus H^1_{Iw, sing}(K_{\infty,\ell},T_f/p^nT_f)\]
together with the projections $v_\ell$ and $\partial_\ell$ onto the first and second factor, respectively. Combined with the localisation map $loc_\ell: H^1_{Iw}(K_\infty, T_{f}/p^nT_f)\to  H^1_{Iw}(K_{\infty,\ell},T_f/p^nT_f) $, their first explicit reciprocity law then reads as follows:

\begin{theorem}[{\cite[Thm.\ 4.1]{BD05}}]
  If $\ell$ denotes an $n$-admissible prime, then
  \[\partial_\ell(loc_\ell(\kappa(\ell)))=\mathcal{L}_f \mod p^n\]
  holds in $H^1_{Iw, sing}(K_{\infty,\ell},T_f/p^nT_f)\cong \Lambda/p^n\Lambda $ (up to multiplication by elements of $\mathbb{Z}^\times_p$ and $G_\infty$).
\end{theorem}

Instead of stating also their second explicit reciprocity law, Thm.\ 4.2 in (loc.\ cit.), we just state the following consequence from both, which is Cor.\ 4.3 in (loc.\ cit.):

\begin{corollary}
For all pairs of $n$-admissible primes $\ell_1,\ell_2$ attached to $f,$ the equality
\[v_{\ell_1}(loc_{\ell_1}(\kappa(\ell_2))=v_{\ell_2}(loc_{\ell_2}(\kappa(\ell_1))\]
holds in $H^1_{Iw, fin}(K_{\infty,\ell},T_f/p^nT_f)\cong \Lambda/p^n\Lambda $ (up to multiplication by elements of $\mathbb{Z}^\times_p$ and $G_\infty$).
\end{corollary}

Note that in contrast to most other generalisations  these level-raising reciprocity laws do have some symmetry concerning the primes $\ell_1$ and $\ell_2$,  which reminds one of the original quadratic reciprocity  law.

Moreover, we wonder whether this usage of ''explicit reciprocity laws'' at the beginning of section 4.1 in (loc.\ cit.)

\begin{quote}\it
Both theorems are instances
of explicit reciprocity laws relating these explicit cohomology classes
to special values of L-functions...
\end{quote}

 might be important for the historical question of how the name ''explicit reciprocity law'' got attached to statements like those in subsection \ref{sec:l=p}: once the terminology was established for theorems describing localisations of Euler systems at $\ell \ne p,$ the jump to also using the term in the $\ell = p$ case seems quite natural.

\section{$\varepsilon$-isomorphisms}\label{sec:eps}

Perrin Riou's explicit reciprocity law has an immediate impact towards the  $\varepsilon$-isomorphisms  \`{a} la Kato \cite{kato-lnm}, \cite{kato-lnmII} or Fukaya-Kato \cite{fukaya-kato}. For a survey on these local Tamagawa Number Conjectures, which can be considered as a form of local Iwasawa main conjectures (while the global TNC corresponds to global Iwasawa main conjectures), see \cite{ven-BSD}, \cite{ven-eps}. Here we just want to point out that the formulae of type \ref{thm:BF} are used by Benois and Berger \cite{benois-berger} (strictly speaking they need and use the formulas from \cite{benois} providing information on the integral level), the formulae of type \ref{thm:SV} in \cite{LVZ}, the formulae \ref{thm:nakamura} by Nakamura \cite{nak} (in the cyclotomic case) and  in \cite{MSVW} (in the Lubin-Tate setting) in order to show different cases of this conjecture. The proof of Kato's rank one case \cite{ven-eps} is very much based on the map $Col$ from section \ref{sec:p-adicL}.



\bibliographystyle{amsplain}
\bibliography{../bib/Xbib}

\end{document}